\documentclass[12pt]{amsart}

\advance\textheight\topskip
\usepackage{amsfonts,amsthm,amssymb}
\allowdisplaybreaks
\usepackage[mathscr]{eucal}

\usepackage[dvips,bookmarks=false,hyperfigures=false,a4paper,
  colorlinks=true, linkcolor=black, citecolor=black]{hyperref}

\usepackage[OT2,T1]{fontenc}
\newcommand\cyr{%
  \renewcommand\rmdefault{wncyr}%
  \renewcommand\sfdefault{wncyss}%
  \renewcommand\encodingdefault{OT2}%
  \normalfont
  \selectfont}
\DeclareTextFontCommand{\textcyr}{\cyr}

\usepackage{times}
\usepackage{mathptmx}

\def\XXint#1#2#3{{%
     \setbox0=\hbox{$#1{#2#3}{\int}$}
     \vcenter{\hbox{$#2#3$}}\kern-.5\wd0}}

\newcommand{\W}{\mathsf{W}}
\newcommand{\Srstau}{\mathsf{S}_{rs}(\tau)}
\newcommand{\Urstau}{\mathsf{U}_{rs}(\tau)}

\newcommand{\hA}{\mathscr{A}}
\newcommand{\hB}{\mathscr{B}}
\newcommand{\hC}{\mathscr{C}}
\newcommand{\sfBfirst}{\mathsf{B}^{{\sim}}}
\newcommand{\sfB}{\mathsf{B}}
\newcommand{\sfA}{\mathsf{A}}

\newcommand{\abar}{\overline{a}}
\newcommand{\kbar}{\overline{k}}
\newcommand{\rbar}{\overline{r}}
\newcommand{\jbar}{\overline{j}}
\newcommand{\mbar}{\overline{m}}
\newcommand{\nbar}{\overline{n}}
\newcommand{\sbar}{\overline{s}}

\newcommand{\K}{\kern1pt\mathscr{K}}
\newcommand{\KK}{\mathsf{K}^{\phantom{y}}}
\newcommand{\KKp}{\mathsf{K}'}
\newcommand{\KKK}{\mathbb{K}}
\newcommand{\KKab}[1]{\kern1pt\mathsf{K}^{(#1)}}

\newcommand{\EX}[1]{e^{#1}_{\phantom{H}}}
\newcommand{\ex}[1]{e^{#1}_{\phantom{H}}}

\newcommand{\minorplus}{\pi}
\newcommand{\skewminor}{\varpi}
\newcommand{\tauskewminor}{\varsigma}
\newcommand{\skewmajor}{\varphi}

\newcommand{\SLiiZ}{SL(2,\oZ)}

\newcommand{\Rmin}{\mathscr{R}^{\mathrm{int}}_{p-1}}
\newcommand{\Rpi}{\mathscr{R}_{p+1}}

\def\Charge#1{*{\scriptstyle[#1]}\ar@{}[0,0];[0,0]+<0pt,3pt>;}
\def\Bharge#1{*{\scriptstyle\boldsymbol{[#1]}}\ar@{}[0,0];[0,0]+<0pt,3pt>;}

\def\verma{\bullet\rule[-2pt]{0pt}{9pt}
  \ar@{-}@*{[|(2.5)]}[0,0]+<0pt,0pt>;[0,0]-<10pt,0pt>
  \ar@{-}@*{[|(2.5)]}[0,0]-<1pt,-1pt>;[0,0]+<7pt,-7pt>
  \ar@{}[0,0]+<0pt,0pt>;}

\def\Cverma{\circ\rule[-2pt]{0pt}{9pt}
  \ar@{-}@*{[|(2.5)]}[0,0]+<-3pt,0pt>;[0,0]-<10pt,0pt>
  \ar@{-}@*{[|(2.5)]}[0,0]+<2pt,-2pt>;[0,0]+<6pt,-7pt>
  \ar@{}[0,0]+<0pt,0pt>;}

\def\iverma{\bullet\rule[-2pt]{0pt}{9pt}
  \ar@{-}@*{[|(2.5)]}[]+<0pt,0pt>;[0,0]+<10pt,0pt>
  \ar@{-}@*{[|(2.5)]}[]+<1pt,1pt>;[0,0]+<-7pt,-7pt>
  \ar@{}[0,0]+<0pt,0pt>;}

\def\Civerma{\circ\rule[-2pt]{0pt}{9pt}
  \ar@{-}@*{[|(2.5)]}[]+<3pt,0pt>;[0,0]+<10pt,0pt>
  \ar@{-}@*{[|(2.5)]}[]-<2pt,2pt>;[0,0]+<-6pt,-7pt>
  \ar@{}[0,0]+<0pt,0pt>;}

\def\nothing{*{\mbox{}}\ar@{};[0,0];}

\def\Bullet{*{\rule[-1pt]{0pt}{7pt}\bullet}\ar@{};[0,0];}
\def\Circ{*{\rule[-1pt]{0pt}{7pt}\circ}\ar@{};[0,0];}

\newcommand{\tensor}{\mathbin{\otimes}}

\newcommand{\rep}{\mathscr}

\newcommand{\dd}{\partial}

\newcommand{\bref}[1]{\textbf{\textup{\ref{#1}}}}

\renewcommand{\geq}{\geqslant}
\renewcommand{\leq}{\leqslant}

\newcommand{\hGL}[1]{\widehat{g\ell}(#1)}
\newcommand{\hSSL}[2]{\widehat{s\ell}(#1|#2)}
\newcommand{\hSL}[1]{\widehat{s\ell}(#1)}
\newcommand{\SL}[1]{s\ell(#1)}

\newcommand{\hD}{\widehat{D}(2|1;\alpha)}

\newcommand{\oC}{\mathbb{C}}
\newcommand{\oN}{\mathbb{N}}
\newcommand{\oR}{\mathbb{R}}
\newcommand{\oZ}{\mathbb{Z}}

\newcommand{\simtimesr}{%
  \mathrel{{\times}\kern-2.6pt\raisebox{1.2pt}{\mbox{\tiny $|$}}}}
\newcommand{\simtimesl}{%
  \mathrel{\raisebox{1.2pt}{\mbox{\tiny $|$}}\kern-2.6pt{\times}}}

\newcommand{\mL}{\rep{L}}  %<----- a \hSL2 module
\newcommand{\charSL}[2]{\chi_{{\vphantom{h}#2}}^{\vphantom{y}#1}}
\newcommand{\charW}[2]{\chi_{{\vphantom{h}#2}}^{\vphantom{y}#1}}

\newcommand{\mfrac}[2]{\raisebox{.8pt}{\mbox{\small$\displaystyle\frac{#1}{#2}$}}}
\newcommand{\ffrac}[2]{\raisebox{.5pt}{\mbox{\footnotesize$\displaystyle\frac{#1}{#2}$}}}
\newcommand{\fffrac}[2]{\raisebox{.9pt}{\mbox{\scriptsize$\displaystyle\frac{#1}{#2}$}}}
\newcommand{\half}{%
  \mathchoice{\ffrac{1}{2}}{\frac{1}{2}}{\frac{1}{2}}{\frac{1}{2}}}

\numberwithin{equation}{section}
\makeatletter
\@addtoreset{equation}{section}
\@addtoreset{subsubsection}{section}

\@addtoreset{figure}{section}

\def\@secnumfont{\bfseries}
\def\subsubsection{\@startsection{subsubsection}{3}%
  \z@{.5\linespacing\@plus.7\linespacing}{-.5em}%
  {\normalfont\bfseries}}
\def\paragraph{\@startsection{paragraph}{4}%
  \z@\z@{-\fontdimen2\font}%
  \normalfont\bfseries}
\def\subparagraph{\@startsection{subparagraph}{5}%
  \z@\z@{-\fontdimen2\font}%
  \normalfont\bfseries}

\makeatother

\swapnumbers

\newtheorem{thm}[subsubsection]{Theorem}
\newtheorem{Lemma}[subsection]{Lemma}
\newtheorem{lemma}[subsubsection]{Lemma}

\newtheorem{cor}[subsubsection]{Corollary}

\theoremstyle{definition}

\newtheorem{rem}[subsubsection]{Remark}

\begin{document}
\hfuzz=2pt

\title[Logarithmic \lowercase{$\hSL2/u(1)$} and Appell functions]{%
  Higher string functions, higher-level Appell functions, and the
  logarithmic $\hSL2_k/u(1)$ CFT model}

\author[Semikhatov]{A.\,M.~Semikhatov}%

\address{\mbox{}\kern-\parindent Lebedev Physics Institute
  \hfill\mbox{}\linebreak \texttt{ams@sci.lebedev.ru}}

\begin{abstract}
  \ We \,generalize \,the \,string \,functions \;$\hC_{n,r}(\tau)$
  \;associated with the coset $\hSL2_k/u(1)$ to higher string
  functions $\hA_{n,r}(\tau)$ and $\hB_{n,r}(\tau)$ associated with
  the coset $\W(k)/u(1)$ of the $W$-algebra of the logarithmically
  extended $\hSL2_k$ conformal field model with positive integer~$k$.
  The higher string functions occur in decomposing $\W(k)$ characters
  with respect to \hbox{level-$k$} theta and Appell functions and
  their derivatives (the characters are neither quasiperiodic nor
  holomorphic, and therefore cannot decompose with respect to only
  theta-functions).  The decomposition coefficients, to be considered
  ``logarithmic parafermionic characters,'' are given by
  $\hA_{n,r}(\tau)$, $\hB_{n,r}(\tau)$, $\hC_{n,r}(\tau)$, and by the
  triplet $\mathscr{W}(p)$-algebra characters of the $(p=k+2,1)$
  logarithmic model.  We study the properties of $\hA_{n,r}$ and
  $\hB_{n,r}$, which nontrivially generalize those of the classic
  string functions~$\hC_{n,r}$, and evaluate the modular group
  representation generated from $\hA_{n,r}(\tau)$ and
  $\hB_{n,r}(\tau)$; its structure inherits some features of modular
  transformations of the higher-level Appell functions and the
  associated transcendental function~$\Phi$.
\end{abstract}

\maketitle
\thispagestyle{empty}

\vspace*{-1.5\baselineskip}

\enlargethispage{1.2\baselineskip}

\begin{footnotesize}\addtolength{\baselineskip}{-6pt}
  \tableofcontents
\end{footnotesize}

\begin{flushright}\footnotesize
  \textcyr{\it Tam dazhe krasovalasp1 neyasnym logarifmom}\\
  \textcyr{\it Abstraktnaya kartina dlya obshchei0 krasoty.}\\
  \textcyr{{Yu.~Vizbor}%%%, {\it``V Arkashinoi0 kvartire.''}
  }
\end{flushright}

\section{Introduction}
The defining feature of logarithmic conformal field
theories~\cite{[Sa],[Gurarie],[GK-rat],[G-alg]}, contrasting them from
rational conformal field theories, is the presence of indecomposable
representations of the chiral algebra.  The interesting representation
theory may be considered the basic reason underlying fascinating
features of logarithmic conformal field models and their links with
several related problems, e.g.,
in~\cite{[F-bits],[FFHST],[LMRS],[CF],[FG],[SS],[FGK],[RS],[W-proof],
  [QS],[GR2],[HLZ],[jF],[S-q]}.  In particular, modular group
representations generated from characters in logarithmic models are of
a different structure than the modular group representations occurring
in rational models (cf.~\cite{[F-9596],[FHST],[FGST]}).

In this paper, we derive a modular group representation of a
``logarithmic'' origin, generated from the \textit{higher string
  functions} (for positive integer $k$ and $1\leq r\leq p=k+2$, with
$n-r\in 2\oZ+1$)
\begin{align*}
  \hA_{n,r}(q)
  &= \ffrac{q^{-\frac{n^2}{4k}}}{\eta(q)^2}\sum_{a\in\oZ}\sum_{j\geq 1}(-1)^{j+1}
  \bigl(a-\ffrac{n}{2k}\bigr)^2
  \Bigl(q^{\half j(j-n) + \frac{(2 a p + r)^2}{4p} + \half j(2a p + r)}
  -(r\mapsto-r)\!\Bigr),\\
  \hB_{n,r}(q)
  &= \ffrac{q^{-\frac{n^2}{4k}}}{\eta(q)^2}\sum_{a\in\oZ}\sum_{j\geq 1}(-1)^{j+1}
  \bigl(a-\ffrac{n}{2k}\bigr)
  \Bigl(q^{\half j(j-n) + \frac{(2 a p + r)^2}{4p} + \half j(2a p + r)}
  -(r\mapsto-r)\!\Bigr),
  \\
  \intertext{which generalize the classic string
    functions~\cite{[KP],[JiMi]}}
  \hC_{n,r}(q)
  &= \ffrac{q^{-\frac{n^2}{4k}}}{\eta(q)^2}
  \sum_{a\in\oZ}\sum_{j\geq 1}(-1)^{j+1}
  \Bigl(q^{\half j(j-n) + \frac{(2 a p + r)^2}{4p} + \half j(2a p + r)}
  -(r\mapsto-r)\!\Bigr)
\end{align*}
in an obvious way.  That $\hA_{n,r}$ and $\hB_{n,r}$ can have
reasonable modular properties is not obvious, however, and these
properties are actually nontrivial.  The most striking feature is that
modular $S$-transformations of $\hA_{n,r}$ and $\hB_{n,r}$ involve the
transcendental function%\enlargethispage{\baselineskip}%
\begin{gather}\label{Phi-def}
  \Phi(\tau,\mu)= -\ffrac{i}{2\sqrt{-i\tau}}
  -\half\int_{\oR}\!dx\,\EX{-\pi x^2}\,
  \mfrac{\sinh\Bigl(\!\pi
    x\sqrt{-i\tau}\bigl(1\!+\!2\fffrac{\mu}{\tau}\bigr)\!\Bigr)}{
    \sinh\!\left(\pi x\sqrt{-i\tau}\right)}
\end{gather}  
introduced previously in studying $\hSSL21$ characters~\cite{[STT]}.
\ Less striking but also interesting is that the modular transforms of
$\hB_{n,r}$ and $\hA_{n,r}$ involve $\Phi$ and its derivative
\textit{times the characters} of the $(p,1)$ logarithmic conformal
field model.  The underlying representation-theory reasons are briefly
as follows.

We recall that the string functions $\hC_{n,r}(q)$ are the
coefficients in the decomposition of integrable $\hSL2_k$ characters
with respect to level-$k$ theta-functions.  Their ``logarithmic''
generalizations $\hB_{n,r}(q)$ and $\hA_{n,r}(q)$ occur similarly in
decomposing the characters of a $W$-algebra $\W(k)$ in a
logarithmically extended minimal $\hSL2_k$ conformal field theory
model~\cite{[S-sl2]}; they are thus associated with a logarithmic
extension of $\hSL2_k/u(1)$.  The modular transformations of
$\hA_{n,r}$ and $\hB_{n,r}$ can then be found in much the same way as
in the well-known case with $\hC_{n,r}$.  Both the technical details
(page~\pageref{sec:technical}) and the result
(page~\pageref{sec:results}) make this undertaking interesting.  But
before describing these, we recall some motivation from logarithmic
conformal field theory (we actually need only the characters and their
modular transformations, and therefore some readers may well skip the
next subsection).\footnote{A general context to which the results in
  this paper relate is that of mock theta-functions.  That this
  particular ``mockery'' of theta functions has reasonable properties
  must be traceable to conformal field theory$/$representation theory
  reasons.}

\subsection{Logarithmic conformal field theory background}
The classic string functions $\hC_{n,r}(q)$ are (modulo normalization)
the characters of the coset $\hSL2_k/u(1)$ model\,---\,the
parafermionic theory that could never complain about lack of attention
since its appearance in~\cite{[FZ]} (e.g., see~\cite{[JM]} and the
references therein, \cite{[DQ],[JNS],[N],[LP]} in particular).  The
higher string functions $\hB_{n,r}(q)$ and $\hA_{n,r}(q)$ are
``logarithmic extensions'' of these characters in that they originate
similarly to the $\hC_{n,r}(q)$ from a logarithmically extended
theory.

Logarithmic conformal field theories differ from rational ones in
several ways, the two major effects being as follows.  First, the
chiral space of states of a logarithmic model is the sum not of all
irreducible representations but of all indecomposable projective
modules (cf.~a~discussion in~\cite{[FGST3],[S-sl2],[RS]}). \ Second,
the chiral algebra itself extends to a larger, typically nonlinear
$W$-algebra.  Such extended algebras can be systematically identified
as maximum local algebras acting in the kernel of the differential in
certain complexes associated with screenings.

Logarithmic conformal models can be systematically defined by choosing
a free-field realization, identifying the screenings that select the
(nonextended, to begin with) chiral algebra as their centralizer,
constructing a complex associated with the screenings, and then taking
the kernel of the differential and the maximum local algebra acting
there~\cite{[FHST],[FGST],[FGST3],[S-sl2]}.  

When the nonextended symmetry is the Virasoro algebra, the chiral
algebra is the triplet $W$-algebra
$\mathscr{W}(p)=W_{2,3\times(2p-1)}$~\cite{[K-first],[GK-loc]} for
$(p,1)$ models or a triplet $W$-algebra~\cite{[FGST3]} with generating
currents of dimension $(2p\,{-}\,1)(2p'\,{-}\,1)$ for $(p,p')$ models.
For $(p,1)$ models, in particular, the ``screening-kernel'' approach
yields a ``semi-explicit'' construction~\cite{[FHST],[FGST]} of the
currents generating the $\mathscr{W}(p)$ algebra (in terms of vertex
operators and screenings; also see~\cite{[AM]}) and a description of
its $2p$ irreducible representations, whence their characters follow
as (see~\cite{[F-9596]} for their first derivation)
\begin{equation}\label{1pcharacters}
  \psi^+_{r}(q)
  = \ffrac{r\theta_{r,p}(q) -
    2\theta'_{r,p}(q)}{p\,\eta(q)},\quad
  \psi^-_{r}(q)
  = \ffrac{r\theta_{p-r,p}(q) + 2\theta'_{p-r,p}(q)}{p\,\eta(q)},
  \qquad 1\leq r\leq p.
\end{equation}

When the nonextended symmetry is $\hSL2_k$ with positive integer $k$,
the currents generating the corresponding extended $W$-algebra $\W(k)$
are of dimension $4p\,{-}\,2$ (and charge $\pm(2p\,{-}\,1)$), $p=k+2$
\cite{[S-sl2]}.  But the ``screening-kernel'' approach suffers from a
mismatch between the number of screenings (two) selecting the $\hSL2$
algebra as their centralizer and the number of free fields (three)
entering the free-field construction of $\hSL2$ (a ``runaway''
direction in the $3$-space of vertex-operator momenta is associated
with the spectral~flow).  These two numbers may be equalized by
passing to a coset over $u(1)$, the coset not of the $\hSL2_k$ algebra
as in the nonlogarithmic case but of the extended algebra $\W(k)$ of
the logarithmic model.  Instead of working out the details of the
resulting ``logarithmic parafermion'' model starting from
representation theory, which seems to be quite a laborious task
(cf.~\cite{[JM]} in the nonlogarithmic case), we work at the level of
characters, and this is how the $\hA$ and $\hB$ functions appear.  The
logarithmically extended parafermion model is, strictly speaking,
presently nonexistent beyond as much as can be deduced from its
proposed characters and the modular group representation generated
from them, derived in what follows.

\subsection{Technical issues}\label{sec:technical}In contrast to the
case with the standard string functions, our starting point is given
by the characters not of (the integrable) $\hSL2_k$-representations
but of representations of the extended $W$-alge\-bra $\W(k)$
constructed in~\cite{[S-sl2]}.  The integrable $\hSL2_k$ characters
are quasiperiodic and holomorphic, but the $\W(k)$ characters
are~neither.  The integrable $\hSL2_k$ characters can therefore be
decomposed with respect to a basis of \hbox{level-$k$} theta
functions%%% (which are quasiperiodic and holomorphic)
, yielding the string functions as the decomposition coefficients, but
the $\W(k)$-characters require a larger basis for decomposition and
hence yield more functions as the coefficients.

\begin{itemize}
\item[--]First, the $\W(k)$ characters are expressed in terms of
  theta-functions $\theta_{r,p}(q,z)$ and their derivatives
  $\theta'_{r,p}(q,z)$ and $\theta''_{r,p}(q,z)$; in the
  decomposition, this leads to the occurrence of $\theta_{n,k}(q,z)$,
  $\theta'_{n,k}(q,z)$, and $\theta''_{n,k}(q,z)$, the coefficients
  being $\hA_{n,r}(q)$, $\hB_{n,r}(q)$, and $\hC_{n,r}(q)$.  For the
  higher string functions, the analogue of the well-known periodicity
  $\hC_{n+2k\ell,r}(q)=\hC_{n,r}(q)$ takes a rather remarkable form:
  shifting $n\to n+2k\ell$ gives rise to additional terms containing
  the triplet $\mathscr{W}(p)$-algebra characters $\psi^{\pm}_{r}(q)$,
  with $p=k+2$.  For example,\footnote{The occurrence of
    $\mathscr{W}(p)$ characters may not be very surprising considering
    that the $\hB_{n,r}$ ``remember'' their origin from the $\W(k)$
    algebra whose Hamiltonian reduction is just the $\mathscr{W}(p)$
    algebra~\cite{[S-sl2]}.}
  \begin{align}\label{B-open}
    \hB_{n+2k,r}(q)
    &= \hB_{n,r}(q)
    + \ffrac{\psi^-_{r}(q)}{\eta(q)}\,q^{-\frac{k}{4}(\frac{n}{k} + 1)^2}
    - \ffrac{\psi^+_{r}(q)}{\eta(q)}\,q^{-\frac{k}{4}(\frac{n}{k} + 2)^2}.\\
  \intertext{Generalizations of the ``reflection'' symmetry
  $\hC_{-n,r}(q)=\hC_{n,r}(q)$ also involve these characters, for
  example,}
  \label{B-minus}
  \smash{\hB_{-n,r}(q)}
  &=\smash[t]{
    -\hB_{n,r}(q)-\ffrac{\psi^+_r(q)}{\eta(q)}\,q^{-\frac{n^2}{4k}}}.
  \end{align}

\item[--]Second, because the $\W(k)$ characters are not holomorphic,
  they cannot decompose with respect to theta functions alone; in
  addition to $\theta_{n,k}(q,z)$ and their derivatives, the
  decomposition involves their meromorphic counterparts, the
  \hbox{level-$k$} Appell functions~\cite{[STT]} (also
  see~\cite{[Polisch],[KW-Appell]})
  \begin{equation*}
    \K_k(q^2,z,y)=
    \sum_{m\in\oZ}\mfrac{\displaystyle
      q^{m^2 k} z^{m k}}{\displaystyle
      1 - z\,y\, q^{2m}}.
  \end{equation*}
  Under modular $S$-transformations, they behave as
  \begin{multline}\label{Knew-Stransf-fin}
    \K_k(-\ffrac{1}{\tau},\ffrac{\nu}{\tau},\ffrac{\mu}{\tau})
    {}= \tau \smash{\EX{i\pi k\frac{\nu^2 - \mu^2}{\tau}}}\,
    \K_k(\tau,\nu,\mu)
    \\*
    {}+ \tau \sum_{n=0}^{k-1}
    \EX{i\pi\frac{k}{\tau}(\nu\!+\!\frac{n}{k}\tau)^2}
    \Phi(k\tau,k\mu\!-\!n\tau) \vartheta(k\tau,
    k\nu\!+\!n \tau),
  \end{multline}
  which is the origin of the $\Phi$ function.\footnote{That the
    level-$k$ Appell functions, which were introduced and studied
    in~\cite{[STT]} motivated by their occurrence in some characters
    of the affine Lie superalgebra $\hSSL21$, make their appearance as
    ``decomposition basis'' elements in the $\hSL2/u(1)$ context may
    of course be attributed to the identification (in the supposedly
    rational case at least, see, e.g.,~\cite{[BFST]})
    \begin{equation*}
      \fffrac{\hSL2_k}{u(1)}=\fffrac{\hSSL21_{k'}}{\hGL2_{k'}},
      \qquad
      (k+1)(k'+1)=1.
    \end{equation*}%
  } \ In the decomposition of $\W(k)$ characters, the coefficients at
  the Appell functions are just the $\mathscr{W}(p)$
  characters~$\psi^\pm_r(\tau)$.
\end{itemize}

To summarize, the $\W(k)$ characters, as functions of~$z$, decompose
with respect to level-$k$ theta functions and their first and second
derivatives, and level-$k$ Appell functions and their first
derivatives.  The decomposition coefficients, which are to be
considered the ``log-parafermionic'' characters,
are
\begin{equation}\label{the-set}
  (\psi^\pm_r(\tau), \hC_{n,r}(\tau), \hB_{n,r}(\tau), \hA_{n,r}(\tau))
\end{equation}
with $1\,{\leq}\,r\,{\leq}\,p\,{=}\,k\,{+}\,2$ and
$0\,{\leq}\,n\,{\leq}\,k$, ``modulo'' several relations at the range
boundaries, such as $\hC_{n,p}(\tau)=0$,
$\hB_{0,r}(\tau)=-\frac{\psi^+_r(\tau)}{2\eta(\tau)}$, and
$\hB_{k,r}(\tau)=\frac{\psi^-_r(\tau)}{2\eta(\tau)} -
\frac{\psi^+_r(\tau)}{\eta(\tau)}\, e^{-i\pi\frac{k}{2}\tau}$,
together with $\hC_{n+k,p-r}(\tau) = \hC_{n,r}(\tau)$ (the
$\psi^\pm_r$ actually occur in the combinations
$\psi^\pm_r(\tau)\,e^{-i\pi\frac{n^2}{2k}\tau}/\eta(\tau)$).

\subsection{Results}\label{sec:results} 
The modular group representation generated from the
set~\eqref{the-set} follows from the modular transformations of the
$\W(k)$-algebra characters in~\cite{[S-sl2]} and of the Appell
functions in~\cite{[STT]}.  The simple modular transformation
properties of $\hC_{n,r}(\tau)$ and $\psi^{\pm}_{r}(\tau)$ characters
are of course well known~\cite{[KP],[FHST]}, but $S$-transforms of
$\hB_{n,r}(\tau)$ and $\hA_{n,r}(\tau)$ are new and turn out to
\textit{involve $\psi^{\pm}_{r}(\tau)$ times the $\Phi$
  function}.\enlargethispage{\baselineskip}

\subsubsection{Notation}
We fix an integer $k\geq1$ and set
\begin{equation*}
  p=k+2.
\end{equation*}
The reader is asked to excuse our mixed use of $k$ and
$p$,\pagebreak[3] which sometimes both occur in the same formula; we
frequently use $(-1)^k=(-1)^p$, $k+1=p-1$, and other helpful
identities.  We also use the notation
\begin{equation*}
  \abar=(a\;\mathrm{mod}\;2)\in\{0,1\}
\end{equation*}
for any $a\in\oZ$, and, more generally,
$[a]_{\ell}=(a\;\mathrm{mod}\;\ell)\in\{0,1,\dots,\ell-1\}$.

We resort to the standard abuse by writing $f(\tau,\nu,\mu)$ for
$f(e^{2i\pi\tau},e^{2i\pi\nu},e^{2i\pi\mu})$; it is tacitly assumed
that $q=e^{2i\pi\tau}$ (with $\tau$ in the upper complex half-plane),
$y=e^{2i\pi\mu}$, etc.

\subsubsection{Background}
We first quote the $S$-transform of the triplet $\mathscr{W}(p)$
algebra characters~\cite{[FHST],[FGST]}:
\begin{multline}\label{eq:S-psiplus}
  \psi^+_r(-\ffrac{1}{\tau})=
  \smash[t]{\sqrt{\fffrac{2}{p}}\,
    \ffrac{r}{2p}}
  \bigl((-1)^{r}\psi^+_p(\tau) + \psi^-_p(\tau)\bigr)
  \\[-4pt]
  {}+
  \sqrt{\fffrac{2}{p}}\,
  \sum_{s=1}^{p-1}
  \ffrac{1}{i}\Srstau
  \bigl(\psi^+_s(\tau) + (-1)^{r}\psi^-_s(\tau)\bigr),
\end{multline}
\mbox{}\vspace*{-\baselineskip}
\begin{multline}\label{eq:S-psiminus}
    \psi^-_r(-\ffrac{1}{\tau})=\sqrt{\fffrac{2}{p}}\,
    \fffrac{r}{2p}\bigl((-1)^{p+r}\psi^+_p(\tau) + \psi^-_p(\tau)\bigr)
    \\[-4pt]
    {}+\sqrt{\fffrac{2}{p}}\,
    \sum_{s=1}^{p-1}(-1)^s
    \ffrac{1}{i}\Srstau
    \bigl(\psi^+_s(\tau) + (-1)^{p+r}\psi^-_s(\tau)\bigr),
\end{multline}
where
\begin{equation}\label{Srstau}
  \Srstau=
  \ffrac{i r}{p}\,
  \cos\!\ffrac{\pi r s}{p} + \tau\ffrac{p\!-\!s}{p}\,
  \sin\!\ffrac{\pi r s}{p}.
\end{equation}
A notable feature of logarithmic conformal field theory is the
explicit occurrence of $\tau$ here.  We next recall that the string
functions $\hC_{m,r}(\tau)$ with $\mbar=\overline{r+1}$ \
$S$-transform as~\cite{[KP]}
\begin{equation}\label{modS-C}
  \hC_{m,r}(-\ffrac{1}{\tau})
  =\ffrac{1}{\sqrt{pk}}
  \mathop{\sum_{s=1}^{k+1}\sum_{n=0}^{2 k-1}}\limits_{\nbar=\overline{s+1}}
    e^{i\pi\frac{m\,n}{k}}\sin\!\ffrac{\pi r s}{p}\,
  \hC_{n,s}(\tau).
\end{equation}

The next theorem shows a nontrivial ``merger'' of the above formulas,
additionally incorporating~$\Phi$, in the $S$-transformation of
$\hB_{m,r}(\tau)$.\enlargethispage{1.5\baselineskip}
\begin{thm}\label{thm:B} For %%%%$1\leq m\leq k-1$,
  $1\leq r\leq p$ and $\mbar=\overline{r+1}$, let
  \begin{multline}\label{sfB-def}
    \sfB_{m,r}(\tau)=\hB_{m,r}(\tau)
    -i\ffrac{\psi^+_r(\tau)\!-\!\psi^-_r(\tau)}{
      2\sqrt{-2i k\tau}\,\eta(\tau)}\\
    {}-\ffrac{\psi^+_r(\tau)}{\eta(\tau)}\,\Phi(2k\tau, m\tau)
    +\ffrac{\psi^-_r(\tau)}{\eta(\tau)}\,\Phi(2k\tau, (m\!-\!k)\tau).
  \end{multline}
  Then
  \begin{equation}\label{sfB-properties}
    \begin{aligned}
      &\sfB_{m+2k,r}(\tau)=\sfB_{m,r}(\tau),%%%\quad\ell\in\oZ,
      \\
      &\sfB_{-m,r}(\tau)+\sfB_{m,r}(\tau)
      = 0
    \end{aligned}
  \end{equation}
  and
  \begin{equation*}
    \sfB_{m,r}(-\ffrac{1}{\tau})
    =\ffrac{(-1)^r}{\sqrt{k p}}
    \ffrac{2 r}{p}
    \!\!\!\sum_{\substack{n=1\\ \nbar=\overline{k+1}}}^{k-1}\!\!\!\!
    \sin\!\ffrac{\pi m\,n}{k}\,
    \hB_{-n,p}(\tau)
    -\ffrac{4 i}{\sqrt{k p}}
    \mathop{\sum_{n=1}^{k-1}
      \sum_{s=1}^{p-1}}\limits_{\overline{n+s}=1}
    \sin\!\ffrac{\pi m\,n}{k}\;
    \Srstau
    \hB_{-n,s}(\tau).\pagebreak[3]
  \end{equation*}
\end{thm}

The $\Phi$ functions involved in the $S$-transformation are thus
neatly incorporated in the definition of the
\textit{``$\Phi$-modified'' string functions~\eqref{sfB-def}}, for
which the properties such as~\eqref{B-open} and~\eqref{B-minus} are
``improved,'' to become the respective relations
in~\eqref{sfB-properties}, and the $S$-transform formula takes the
simplest form.  We note that $\sfB_{m,r}(-\frac{1}{\tau})$ (and hence
$\hB_{m,r}(-\frac{1}{\tau})$) with $1\leq m\leq k-1$ are expressed
through $\hB_{n,s}(\tau)$ with $-k+1\leq n\leq-1$; to reexpress the
right-hand side in terms of positive-moded $\hB_{n,s}(\tau)$,
Eqs.~\eqref{B-open}--\eqref{B-minus} must be used; evidently,
expressing the right-hand side in terms of the $\sfB_{m,s}(\tau)$
introduces the $\Phi$ functions.  Iterating the $S$-transformation, in
terms of either $\hB$ or $\sfB$, inevitably leads to accumulating
$\Phi$'s with different arguments, and it is clear that the modular
group relations require that certain such combinations evaluate in
terms of elementary functions (exponentials).  Because $\Phi$ itself
originates from modular transformation~\eqref{Knew-Stransf-fin}, it
satisfies the necessary ``consistency'' conditions, as is detailed
in~\cite{[STT]}; specifically in the string-function context, the
relevant identities are explicitly given in~\bref{sin-Phi-lemma} in
what follows.

The above $S$-transformation may be compared with~\eqref{modS-C},
suggestively rewritten as
\begin{equation*}
  \hC_{m,r}(-\ffrac{1}{\tau})
  = \ffrac{2}{\sqrt{pk}}\mathop{\sum_{n=0}^{k-1}\sum_{s=1}^{p-1}}
  \limits_{\overline{n+s}=1}
  \cos\!\ffrac{\pi m\,n}{k}\,\sin\!\ffrac{\pi r s}{p}\,
  \hC_{-n,s}(\tau).
\end{equation*}
Besides $\sin\frac{\pi m\,n}{k}$ in the theorem replacing
$\cos\frac{\pi m\,n}{k}$ in the above formula (in accordance with the
``odd'' property of $\sfB$ in~\eqref{sfB-properties}), a notable
difference is that $\tau$ explicitly occurs in $\Srstau$, a feature in
common with the $(p,1)$ logarithmic model; but the most essential
increase in complexity in passing to the $\hB$ case is the
incorporation of the $\Phi$ function in~\eqref{sfB-def}.

We also note that $\sfB_{0,r}(q)$ and $\sfB_{k,r}(q)$ defined as
in~\eqref{sfB-def} vanish, which means that $\hB_{k\ell,r}(q)$,
$\ell\in\oZ$, %%%and $\hB_{k,r}(q)$
are in $\oC[\psi^{\pm}_r(\tau)\,e^{-i\pi\frac{n^2}{2
    k}\tau}/\eta(\tau)]$ ($n\in\oZ$%%%, $1\leq s\leq p$
).  The $S$-transform formula in the theorem is therefore consistent
but not informative for $m=0,\;k$.

The proof of the $S$-transform formula in~\bref{thm:B} is the content
of~\bref{sec:S-B}; simple relations \eqref{sfB-properties} are shown
in Appendix~\bref{app:ABC}.\enlargethispage{\baselineskip}

For $\hA_{n,r}$, the counterparts of relations \eqref{B-open} and
\eqref{B-minus} are
\begin{align}\label{A-open}
  \hA_{n + 2 k,r}(q)&=\hA_{n,r}(q)
  -\bigl(1+\ffrac{n}{k}\bigr)
  \ffrac{\psi^-_r(q)}{\eta(q)}\,q^{-\frac{(n+k)^2}{4 k}}
  +\bigl(2+\ffrac{n}{k}\bigr)
  \ffrac{\psi^+_r(q)}{\eta(q)}\,q^{-\frac{(n+2 k)^2}{4 k}}
  \\[-4pt]
  \intertext{and}
  \label{A-minus}
  \hA_{-n,r}(q)&=\hA_{n,r}(q)
  - \smash[t]{\ffrac{n}{k}\ffrac{\psi^+_r(q)}{\eta(q)}\,
    q^{-\frac{n^2}{4 k}}}.
\end{align}
As with the $\hB$, these properties are ``improved'' for
$\Phi$-modified string functions.  We set
\begin{equation}\label{Phi-prime}
  \Phi'(\tau,\mu)=\ffrac{1}{2i\pi}\,\ffrac{\dd}{\dd\mu}\,
  \Phi(\tau,\mu).
\end{equation}
\begin{thm}\label{thm:A}
  For %%%%%$1\leq m\leq k-1$,
  $1\leq r\leq p$ and $\mbar=\overline{r+1}$,
  let
  \begin{equation}\label{sfA-def}
    \sfA_{m,r}(\tau) = \hA_{m,r}(\tau)
    - \ffrac{2\psi^+_{r}(\tau)}{\eta(\tau)}\,
    \Phi'\bigl(2k\tau, m\tau\bigr)
    + \ffrac{2\psi^-_{r}(\tau)}{\eta(\tau)}\,
    \Phi'\bigl(2k\tau, (m - k)\tau\bigr).
  \end{equation}
  Then
  \begin{equation}\label{sfA-properties}
    \begin{aligned}
      &\sfA_{m+2k,r}(\tau)=\sfA_{m,r}(\tau),
      \\
      &\sfA_{-m,r}(\tau)-\sfA_{m,r}(\tau)=0,
    \end{aligned}
  \end{equation}
  and
  \begin{multline*}
    \sfA_{m,r}(-\ffrac{1}{\tau})=\\[-6pt]
    {}=\ffrac{2\tau}{\sqrt{p\,k\!}}
    \Bigl(\!
    (1+(-1)^{k+1})\ffrac{(-1)^r i r\!}{4 p}
    \hA_{0,p}(\tau)  
    +\!\sum_{\substack{s=1\\ \sbar=1}}^{p-1}\!
    \Srstau \hA_{0,s}(\tau)
    - (-1)^r
    \!\!\!\!\sum_{\substack{s=1\\ \sbar=\overline{k+1}}}^{p-1}\!\!\!\!
    \Srstau \hA_{-k,s}(\tau)\!\!\Bigr)
    \\[-4pt]
    +\ffrac{4\tau}{\sqrt{p\,k}}
    \!\!\sum_{\substack{n=1\\ \nbar=\overline{k+1}}}^{k-1}\!\!
    \cos\!\ffrac{\pi m\,n}{k}\;
    \ffrac{(-1)^r i r}{2 p}\,\hA_{-n,p}(\tau)
    +\ffrac{4\tau}{\sqrt{p\,k}}\!
    \mathop{\sum_{n=1}^{k-1}\sum_{s=1}^{p-1}}\limits_{\overline{n+s}=1}
    \cos\!\ffrac{\pi m\,n}{k}\;
    \Srstau \hA_{-n,s}(\tau)
    \\[-4pt]
    +\ffrac{1}{\sqrt{p\,k}}\!
    \mathop{\sum_{n=0}^{k-1}\sum_{s=1}^{p-1}}\limits_{\overline{n+s}=1}
    \cos\!\ffrac{\pi m\,n}{k}\;
    \Urstau\,\hC_{n,s}(\tau),
  \end{multline*}
  where
  \begin{equation*}
    \Urstau=\ffrac{i r (p - 2 s)\tau}{2 p^2}\,
    \cos\!\ffrac{\pi r s}{p}
    +\Bigl(\ffrac{s(s-p)\tau^2}{2 p^2} 
    -\ffrac{\tau}{i\pi p\,k}
    + \ffrac{r^2}{2 p^2}
    \Bigr)
    \sin\!\ffrac{\pi r s}{p}.
  \end{equation*}
\end{thm}
This formula looks more complicated than its ``lower'' analogue
in~\bref{thm:B} for three reasons: $\sfA_{k\ell,r}$, $\ell\in\oZ$, do
not vanish and hence contribute to the transformation; also, the
``$\cos\frac{\pi m n}{k}$'' representation of $\SLiiZ$ is somewhat
bulkier than the ``$\sin\frac{\pi m n}{k}$'' representation (when
$\hA_{0,*}$ are not related to $\hA_{\pm k,*}$); finally, there is an
``admixture'' of the $\hC$ string functions.

The proof of the $S$-transform formula is the content
of~\bref{sec:S-B}; simple relations~\eqref{sfA-properties} are shown
in Appendix~\bref{app:ABC}.

\subsubsection{}We note that the $T$-transformation
$\tau\mapsto\tau+1$ amounts to multiplying $\hA_{m,r}(\tau)$,
$\hB_{m,r}(\tau)$, and $\hC_{m,r}(\tau)$ by $e^{i\pi(\frac{r^2}{2 p} -
  \frac{n^2}{2 k} - \frac{1}{6})}$ and
$\fffrac{\psi^{\pm}_r(\tau)}{\eta(\tau)}$ by $e^{i\pi(\frac{r^2}{2 p}
  - \frac{1}{6})}$.

\subsubsection*{Plan of the paper} We extract the higher string
functions from decomposing the characters of the triplet $W$-algebra
$\W(k)$ of logarithmically extended $\hSL2_k$-models in
Sec.~\bref{sec:decomposition}.  Modular $S$-trans\-for\-mations of the
higher string functions are derived in Sec.~\bref{sec:modular}.
Theta-function conventions are fixed in Appendix~\bref{app:theta}.
The necessary properties of the Appell functions are recalled in
Appendix~\bref{sec:Appell}.  The $\W(k)$-algebra characters are listed
and their modular properties are recalled in
Appendix~\bref{sec:Wchars}.  Some simple properties of the higher
string functions are derived in Appendix~\bref{app:ABC}.

The calculations leading to the results stated above are
straightforward but quite bulky.  Besides, the Appell functions
$\KK\!\!$ and the related $\Phi$ function are integrated into the
derivation, and their properties have a considerable impact on the
``calculation flow,'' with the ``sign'' of the effect dependent on
whether these properties are used timely or untimely.  Essential
simplification (although possibly still far from the ideal) is
achieved by consolidating the relevant $\KK\!\!\!/\Phi$ properties
in~\bref{lemma:basic-S}.

\section{Character decompositions}\label{sec:decomposition} 
In this section, we establish the decomposition, or ``branching,'' of
the $\W(k)$-algebra characters in~\cite{[S-sl2]}.  The method is very
direct and is based on the identity (see~\cite{[SW],[KP]} and the
references therein)
\begin{equation}\label{the-identity}
  \ffrac{1}{q^{\frac{1}{8}}\,\vartheta_{1,1}(q,z)}=
  \ffrac{1}{\eta(q)^3}
  \sum_{m\in\oZ}\sum_{j\geq1}(-1)^{j+1}q^{\half j(j-1) - j m}z^{-m}.
\end{equation}

\subsection{$\smash{\hSL2}$ integrable representation characters} We
first recall the classic result~\cite{[KP],[JiMi]} that the integrable
$\hSL2_k$ characters
\begin{equation}\label{integrable-sl2}
  \charSL{}{r}(q,z)
  = \ffrac{\theta_{r,p}(q,z) - \theta_{-r,p}(q,z)}{
    \Omega(q,z)
  },\quad r=1,\dots,k+1,
\end{equation}
decompose with respect to level-$k$ theta-functions as
\begin{equation}\label{int-decomp}
  \charSL{}{r}(q,z)
  =-\ffrac{1}{\eta(q)}
  \sum_{\substack{n=0\\ \nbar=\overline{r+1}}}^{2 k-1}\hC_{n,r}(q)
  \theta_{n,k}(q,z).
\end{equation}
Theta-function conventions and the definition of $\Omega(q,z)$ are
given in Appendix~\ref{app:theta}.

We next decompose the other $\W(k)$-characters similarly
to~\eqref{int-decomp}.

\subsection{Decomposition of the $\W(k)$-algebra characters}
\subsubsection{The characters}\label{sec:chars}
In the logarithmic $\hSL2_k$ model for each $k=0,1,2,\dots$,
characters of the extended algebra $\W(k)$ were calculated
in~\cite{[S-sl2]}.  The characters $\charW{\pm}{r}(q,z)$ are given
by\enlargethispage{\baselineskip}
\begin{align}\label{chi-plus}
  \charW{+}{r}(q,z) &=\ffrac{1}{\Omega(q,z)}\biggl(
  \ffrac{r^2}{4p^2}\bigl(\theta_{-r,p}(q,z)-\theta_{r,p}(q,z)\bigr)
  \\*
  {}&\qquad\quad{}
  +\ffrac{r}{p^2}\bigl(\theta'_{-r,p}(q,z)+\theta'_{r,p}(q,z)\bigr)
  +\ffrac{1}{p^2}\bigl(\theta''_{-r,p}(q,z)-\theta''_{r,p}(q,z)\bigr)
  \!\!\biggr),
  \notag\\
  \label{chi-minus}
  \charW{-}{r}(q,z) &=\ffrac{1}{\Omega(q,z)}\biggl(\!\!
  \Bigl(\ffrac{r^2}{4p^2} -
  \ffrac{1}{4}\Bigr)\bigl(\theta_{p-r,p}(q,z)-\theta_{p+r,p}(q,z)\bigr)
  \\*
  {}&\qquad{}
  +\ffrac{r}{p^2}\bigl(\theta'_{p-r,p}(q,z)+\theta'_{p+r,p}(q,z)\bigr)
  +\ffrac{1}{p^2}\bigl(\theta''_{p-r,p}(q,z)-\theta''_{p+r,p}(q,z)\bigr)
  \!\!\biggr)
  \notag
\end{align}
for $1\leq r\leq p\,{-}\,1$ and
\begin{equation*}
  \charW{+}{p}(q,z)
  = \ffrac{2 \theta'_{p,p}(q,z)}{p\,\Omega(q,z)},
  \qquad
  \charW{-}{p}(q,z)
  =\ffrac{2 \theta'_{0,p}(q,z)}{p\,\Omega(q,z)}.
\end{equation*}

Under the spectral flow (see~\bref{sec:sf}), the $\charW{\pm}{r}(q,z)$
further generate~$\omega^{\pm}_r(q,z)$ given by~\cite{[S-sl2]}
\begin{align*}
  \omega^+_r(q,z)
  &=\ffrac{1}{\Omega(q,z)}\Bigl(
  %%%+
  \ffrac{r}{2p}\bigl(\theta_{r,p}(q,z)+\theta_{-r,p}(q,z)\bigr)
  -
  \ffrac{1}{p}\bigl(\theta'_{r,p}(q,z)-\theta'_{-r,p}(q,z)\bigr)
  \!\Bigr),\\
  \omega^-_r(q,z)
  &=\ffrac{1}{\Omega(q,z)}\Bigl(
  %%%+
  \ffrac{r}{2p}\bigl(\theta_{p-r,p}(q,z)+\theta_{r-p,p}(q,z)\bigr)
  -
  \ffrac{1}{p}\bigl(\theta'_{r-p,p}(q,z)-\theta'_{p-r,p}(q,z)\bigr)
  \!\Bigr)
  \\
  \intertext{for $ 1\leq r\leq p\!-\!1$, and}
  \omega^+_p(q,z)&=
  \ffrac{\theta_{p,p}(q,z)}{\Omega(q,z)},\qquad
  \omega^-_p(q,z)=
  \ffrac{\theta_{0,p}(q,z)}{\Omega(q,z)}.
\end{align*}

The characters decompose with respect to level-$k$ theta and Appell
functions and their derivatives.  We set
\begin{equation*}
  \KKp_{\alpha,k}(q,x,y)
  =\bigl(x\ffrac{\dd}{\dd x} - y\ffrac{\dd}{\dd y}\bigr)
  \KK_{\alpha,k}(q,x,y),
\end{equation*}
where the functions $\KK_{\alpha,k}(q,x,y)$ defined in~\eqref{eq:KK}.
\begin{lemma}\label{chi-lemma}
  As functions of $z$, the $\W(k)$ characters $\charW{\pm}{r}(q,z)$
  decompose with respect to level-$k$ theta and Appell functions and
  their derivatives as
  \begin{multline*}
    \charW{+}{r}(q,z)
    =
    \ffrac{1}{\eta(q)}
    \!\!\sum_{\substack{n=0\\ \nbar=\overline{r+1}}}^{2 k-1}\!\!
    \Bigl(
    \hA_{n,r}(q)\theta_{n,k}(q,z)
    + \ffrac{2}{k}\,\hB_{n,r}(q)\theta'_{n,k}(q,z)
    + \ffrac{1}{k^2}\,\hC_{n,r}(q)\theta''_{n,k}(q,z)
    \!\Bigr)
    \\*
    {}-\ffrac{2\psi^+_{r}(q)}{\eta(q)^2}\,
     q^{-k}\Bigl(\ffrac{1}{k}\,\KKp_{\overline{r+1},k}(q,z,q^{-2})
    - \KK_{\overline{r+1},k}(q,z,q^{-2})\Bigr)
    \\
    {}+ \ffrac{\psi^-_{r}(q)}{\eta(q)^2}\,
    q^{-\frac{k}{4}}\Bigl(\ffrac{2}{k}\,\KKp_{\overline{r+1},k}(q,z,q^{-1})
    - \KK_{\overline{r+1},k}(q,z,q^{-1})
    \!\Bigr)
  \end{multline*}
  and
  \begin{multline*}
    \charW{-}{r}(q,z)
    =-\ffrac{1}{\eta(q)}
    \smash[t]{\!\!\sum_{\substack{n=0\\ \nbar=\overline{r+1}}}^{2 k-1}\!\!}
    \Bigl(\!\!
    \bigl(\hA_{n,r}(q)-\ffrac{1}{4}\,\hC_{n,r}(q)\bigr)
    \theta_{n+k,k}(q,z)
    \\[-9pt]
    \shoveright{{}
      +\ffrac{2}{k}\,\hB_{n,r}(q)\theta'_{n+k,k}(q,z)
      +\ffrac{1}{k^2}\,\hC_{n,r}(q)\theta''_{n+k,k}(q,z)
      \!\Bigr)}
    \\
    {}+\ffrac{\psi^+_{r}(q)}{\eta(q)^2}\,
    q^{-\frac{k}{4}}\Bigl(\ffrac{2}{k}\,\KKp_{k+\overline{r+1},k}(q,z,q^{-1})
    - \KK_{k+\overline{r+1},k}(q,z,q^{-1})
    \!\Bigr)\\
    {}-\ffrac{\psi^-_{r}(q)}{\eta(q)^2}
    \,\ffrac{2}{k}\,\KKp_{k+\overline{r+1},k}(q,z,1).
  \end{multline*}
%%   \textup{(}where $\psi^{\pm}_{r}(q)$ are the $\mathscr{W}(p)$
%%   characters in~\eqref{1pcharacters}\textup{)}.
\end{lemma}

A simple corollary follows if we use~\eqref{minus-twisted} to evaluate
$\omega^+_r(q,z)=-\charW{-}{r;1}(q,z) - \charW{+}{r}(q,z)$ and
$\omega^-_r(q,z) =-\charW{+}{r;1}(q,z) - \charW{-}{r}(q,z)
-\half\charSL{}{p-r}(q,z)$ with the above decompositions of
$\charW{\pm}{r}$ and~$\charSL{}{r}$.

\begin{cor}\label{lemma:minor}There are the decompositions
  \begin{multline*}
    \omega^+_r(q,z)
    =\ffrac{1}{\eta(q)}
    \!\!\sum_{\substack{n=0\\ \nbar=\overline{r+1}}}^{2 k-1}\!\!
    \Bigl(
    \hB_{n,r}(q)\theta_{n,k}(q,z)
    +\ffrac{1}{k}\,\hC_{n,r}(q)\theta'_{n,k}(q,z)
    \!\!\Bigr)
    \\*[-7pt]
    {}-\ffrac{\psi^+_{r}(q)}{\eta(q)^2}\,
    q^{-k}\,\KK_{\overline{r+1},k}(q,z,q^{-2})
    +
    \ffrac{\psi^-_{r}(q)}{\eta(q)^2}\,
    q^{-\frac{k}{4}}\,\KK_{\overline{r+1},k}(q,z,q^{-1})
  \end{multline*}
  and
  \begin{multline*}%%%\label{minor-minus}
    \omega^-_r(q,z)
    {}=
    -\ffrac{1}{\eta(q)}
    \!\!\smash[t]{\sum_{\substack{n=0\\ \nbar=\overline{r+1}}}^{2 k-1}}\!\!
    \Bigl(
    \hB_{n,r}(q)\theta_{n+k,k}(q,z)
    + \ffrac{1}{k}\,
    \hC_{n,r}(q)\theta'_{n+k,k}(q,z)
    \!\!\Bigr)
    \\*[-7pt]
    {}+\ffrac{\psi^+_{r}(q)}{\eta(q)^2}
    \,q^{-\frac{k}{4}}\,\KK_{k+\overline{r+1},k}(q,z,q^{-1})
    - \ffrac{\psi^-_{r}(q)}{\eta(q)^2}
    \,\KK_{k+\overline{r+1},k}(q,z,1).
  \end{multline*}
\end{cor}
%% The rest of this section is the proof of the two decompositions
%% in~\bref{chi-lemma}.

\subsection{Proof of~\bref{chi-lemma}}\label{app:proof-chi}
We derive the decomposition formula for $\charW{-}{r}(q,z)$
in~\bref{app:proof-minus} and the formula for $\charW{+}{r}(q,z)$
in~\bref{app:proof-plus}.

\subsubsection{$\protect\charW{-}{r}(q,z)$}\label{app:proof-minus}
We write the $\charW{-}{r}$ character in~\eqref{chi-minus} as
\begin{equation*}
  \charW{-}{r}(q,z)
  =\ffrac{1}{q^{\frac{1}{8}}\,\vartheta_{1,1}(q,z)}
  \sum_{a\in\oZ+\half}
  (a^2 - \ffrac{1}{4})q^{p(\frac{r}{2p}+a)^2}\bigl(z^{-\frac{r+1}{2}-ap}
  -z^{\frac{r-1}{2}+ap}\bigr).
\end{equation*}
Using identity~\eqref{the-identity}, we calculate
\begin{multline*}
  \charW{-}{r}(q,z)=\\
  {}=\mfrac{q^{\frac{r^2}{4p}}}{\eta(q)^3}\!\!
  \sum_{a\in\oZ+\half}\sum_{m\in\oZ}\sum_{j\geq1}
  (-1)^{j+1}\bigl(a^2-\ffrac{1}{4}\bigr)
  q^{\half j(j-1) - j m + p a^2 + r a}
  z^{-\frac{r}{2}-a p - m - \half}
   - (r\mapsto-r).
\end{multline*}
We now shift the summation variable as $a\mapsto a-\half$ and then
pass from summation over $m$ to summation over $n=2 m + r + 2a p - 1$,
which (with integer-valued $a$) ranges over $2\oZ+\overline{r+1}$.
Shifting $j\mapsto j+1$ then yields
\begin{multline*}
  \charW{-}{r}(q,z)=\\
  {}=\ffrac{q^{\frac{r^2}{4p}}}{\eta(q)^3}\!\!\!
  \sum_{n\in 2\oZ+\overline{r+1}}\!\!\!\!
  q^{\frac{k}{4}-\frac{n}{2}} z^{\frac{k}{2}-\frac{n}{2}}
  \sum_{a\in\oZ}\sum_{j\geq0}\!
  (-1)^jq^{\half j(j-n)}(a^2\!-\!a)
  q^{j(a p + \frac{r}{2}) + p a^2 + r a}
  - (r\mapsto-r).
\end{multline*}

Next, the elementary identity
\begin{equation}\label{zero-dentity}
  \sum_{j\in\oZ}(-1)^j q^{\half j(j-1) + j n} = 0,\quad n\in\oZ,
\end{equation}
and the antisymmetry of the entire expression for $\charW{-}{r}$ under
$r\mapsto-r$ allow us to conclude that
\begin{equation*}
  \charW{-}{r}(q,z)
  =-\ffrac{q^{\frac{r^2}{4p}}}{\eta(q)^3}\!\!
  \sum_{n\in 2\oZ+\overline{r+1}}
  q^{\frac{k}{4} + \frac{n}{2}}\,z^{\frac{k}{2} + \frac{n}{2}}
  (A_{n,r}(q) + B_{n,r}(q))
\end{equation*}
with $A_{n,r}(q)$ and $B_{n,r}(q)$ defined in~\eqref{the-A}
and~\eqref{the-B}.  The formulas in~\bref{app:theF} for
$A_{n+2k\ell,r}$ and $B_{n+2k\ell,r}$ then yield
\begin{multline}\label{Psi-part}
  \charW{-}{r}(q,z)
%%   =-\ffrac{1}{\eta(q)}
%%   \!\smash[t]{\sum_{b=0}^{k-1}}\Bigl(\!
%%   \bigl(\hA_{2b+\overline{r+1},r}(q)
%%   -\ffrac{1}{4}\,\hC_{2b+\overline{r+1},r}(q)\bigr)
%%   \theta_{2b+\overline{r+1}+k,k}(q,z)\\*
%%   \shoveright{{}
%%     +\ffrac{2}{k}\,\hB_{2b+\overline{r+1},r}(q)
%%     \theta'_{2b+\overline{r+1}+k,k}(q,z)
%%     +\ffrac{1}{k^2}\,\hC_{2b+\overline{r+1},r}(q)
%%     \theta''_{2b+\overline{r+1}+k,k}(q,z)
%%     \!\Bigr)\;}
%%   \\
%%   {}+\ffrac{1}{\eta(q)^2}\!
%%   \sum_{b=0}^{k-1}
%%   \Bigl(\sum_{j\geq1}\sum_{\ell\geq j}-\sum_{j\leq0}\sum_{\ell\leq j-1}\Bigr)
%%   q^{\frac{k}{4}+b+\half(\overline{r+1}) + k\ell^2 + k\ell + (2b+\overline{r+1})\ell}\,
%%   z^{\frac{k}{2} + b + \half(\overline{r+1}) + k\ell}\times{}\\
%%   {}\times
%%   \Bigl(\!(2\ell-2j+1)q^{-k j^2 - j(2b+\overline{r+1})}\psi^+_{r}(q)\\
%%   {}-(2\ell-2j+2)
%%   q^{-\frac{k}{4}(2j-1)^2 - (j-\half)(2b+\overline{r+1})}\psi^-_{r}(q)
%%   \!\Bigr),
  =-\ffrac{1}{\eta(q)}
  \!\smash[t]{%
    \sum_{\substack{n=0\\ \nbar=\overline{r+1}}}^{2 k-1}}\Bigl(\!
  \bigl(\hA_{n,r}(q)
  -\ffrac{1}{4}\,\hC_{n,r}(q)\bigr)
  \theta_{n+k,k}(q,z)\\[-6pt]
  \shoveright{{}
    +\ffrac{2}{k}\,\hB_{n,r}(q)
    \theta'_{n+k,k}(q,z)
    +\ffrac{1}{k^2}\,\hC_{n,r}(q)
    \theta''_{n+k,k}(q,z)
    \!\Bigr)\;}
  \\
  \shoveleft{\quad
    {}+\ffrac{1}{\eta(q)^2}\!
    \sum_{\substack{n=0\\ \nbar=\overline{r+1}}}^{2 k-1}
    \Bigl(\sum_{j\geq1}\sum_{\ell\geq j}-\sum_{j\leq0}\sum_{\ell\leq j-1}\Bigr)
    q^{\frac{k}{4} + \frac{n}{2}  + k\ell^2 + k\ell + (n)\ell}\,
    z^{\frac{k}{2} + \frac{n}{2} + k\ell}\times{}}
  \\
  {}\times
  \Bigl(\!(2\ell-2j+1)q^{-k j^2 - j(n)}\psi^+_{r}(q)
  -(2\ell-2j+2)
  q^{-\frac{k}{4}(2j-1)^2 - (j-\half)(n)}\psi^-_{r}(q)
  \!\Bigr),
\end{multline}
where
\begin{align}
  \hC_{n,r}(q)&=\smash[t]{\ffrac{q^{\frac{r^2}{4p}-\frac{n^2}{4k}}}{\eta(q)^2}}
  \,C_{n,r}(q),\notag\\
  \label{ABC-help}
  \hB_{n,r}(q)&=\ffrac{q^{\frac{r^2}{4p}-\frac{n^2}{4k}}}{\eta(q)^2}\bigl(
  B_{n,r}(q)-\ffrac{n}{2 k}\,C_{n,r}(q)\bigr),\\
  \hA_{n,r}(q)&=\ffrac{q^{\frac{r^2}{4p}-\frac{n^2}{4k}}}{\eta(q)^2}\bigl(
  A_{n,r}(q)-\ffrac{n}{k}\,B_{n,r}(q)
  +\ffrac{n^2}{4 k^2}\,C_{n,r}(q)\bigr),\notag
\end{align}

In the ``$\psi$-part'' of~\eqref{Psi-part}, we make the shift
$\ell\mapsto\ell+j$, which produces the sums
\begin{equation*}
  \sum_{j\geq1}\sum_{\ell\geq0}-\sum_{j\leq0}\sum_{\ell\leq-1}
  =\sum_{\ell\geq0}\sum_{j\geq1}-\sum_{\ell\leq-1}\sum_{j\leq0},
\end{equation*}
and evaluate the resulting $j$-sums.  Examination shows that under the
condition
\begin{equation}\label{qz1}
  |q|<|z|<1,
\end{equation}
all of the $j$-sums are of the form $\sum_j\xi^j$ with $|\xi|<1$,
summed over positive (nonnegative) $j$.  For each $\ell\leq-1$, for
instance, the coefficient at $\psi^-_{r}(q)$ involves the sums
\begin{equation*}
  \sum_{j\leq0}(q^{2k(\ell+1)}\,z^k)^j
  = \sum_{j\geq0}(q^{2k(-\ell-1)}\,z^{-k})^j,
\end{equation*}
where $|q^{2k(-\ell-1)}\,z^{-k}|<1$ for any $\ell\leq-2$.  This
estimate does not hold in the sole case $\ell=-1$, but the divergent
sum $\sum_{j\geq0}z^{-k j}$ does not actually occur because of the
factor $(2\ell+2)$ in front of $\psi^-_{r}(q)$ (after the shift
$\ell\mapsto\ell+j$ in the ``$\psi$-part'' in~\eqref{Psi-part}).  The
result~is
\begin{multline*}
  \charW{-}{r}(q,z) =-\ffrac{1}{\eta(q)}\!\!
  \sum_{\substack{n=0\\ \nbar=\overline{r+1}}}^{2 k-1}\Bigl(\!
  \bigl(\hA_{n,r}(q)-\ffrac{1}{4}\,\hC_{n,r}(q)\bigr)
  \theta_{n+k,k}(q,z)
  \\[-6pt]
  \shoveright{{}
    +\ffrac{2}{k}\,\hB_{n,r}(q)\theta'_{n+k,k}(q,z)
    +\ffrac{1}{k^2}\,\hC_{n,r}(q)\theta''_{n+k,k}(q,z)
    \!\!\Bigr)\;}
  \\
  \shoveleft{\phantom{\charW{-}{r}(q,z) ={}}
    +\ffrac{1}{\eta(q)^2}\!\sum_{\ell\in\oZ}
    q^{\frac{k}{4}+\half(\overline{r+1}) + k\ell^2 + k\ell
      + (\overline{r+1})\ell}q^{2(\ell+1)k}
    \,z^{\frac{k}{2}+\half(\overline{r+1}) + k\ell+k}\times{}}
  \\[-6pt]
  {}\times
  \Bigl(
  \mfrac{(2\ell+1)q^{-k}}{1-z\,q^{2\ell+1}}\,\psi^+_{r}(q)
  -\mfrac{(2\ell+2)q^{-\frac{k}{4}+\half(\overline{r+1})}}{
    1-z\,q^{2(\ell+1)}}\,\psi^-_{r}(q)
  \!\!\Bigr).
\end{multline*}
After simple rearrangements, we obtain the formula in the theorem.

\subsubsection{$\protect\charW{+}{r}(q,z)$}\label{app:proof-plus}
We write the $\charW{+}{r}$ character as
\begin{equation*}
  \charW{+}{r}(q,z)
  =\ffrac{1}{q^{\frac{1}{8}}\,\vartheta_{1,1}(q,z)} \sum_{a\in\oZ}
  a^2q^{p(\frac{r}{2p}+a)^2}\bigl(z^{-\frac{r+1}{2}-ap}
  -z^{\frac{r-1}{2}+ap}\bigr).\pagebreak[3]
\end{equation*}
Using~\eqref{the-identity} again, we find
\begin{multline*}%%%\label{char-minus}
  \charW{+}{r}(q,z)
%%   =\\
  {}=
  \ffrac{1}{\eta(q)^3}
  \!\sum_{m\in\oZ}\!z^{-m-\half}
  \sum_{j\geq1}\!(-1)^{j+1}q^{\half j(j-1) - j m}
  \sum_{a\in\oZ}
  \sum_{\sigma=\pm1}\!\!\!
  \sigma a^2q^{p(\frac{r}{2p}+a)^2} z^{-\sigma\frac{r}{2}-\sigma ap}\\
  =\mfrac{q^{\frac{r^2}{4p}}}{\eta(q)^3}
  \sum_{n\in 2\oZ+\overline{r+1}}\!\!\!
  z^{-\frac{n}{2}}
  A_{n,r}(q),
\end{multline*}
with $A_{n,r}(q)$ defined in~\eqref{the-A}.

Next, identity~\eqref{zero-dentity} shows that
\begin{multline*}
  A_{n,r}-A_{-n,r}=\\
  =\sum_{a\in\oZ}a^2
  \Bigl(\sum_{j\geq1}+\sum_{j\leq-1}\Bigr) 
   (-1)^{j+1} q^{\half j(j-n) + r a + p a^2}
  \Bigl(q^{\half j(2a p + r)} - q^{-\half j(2 a p + r)}\Bigr)=0,
\end{multline*}
and therefore
\begin{equation*}
  \sum_{n\in 2\oZ+\overline{r+1}}\!\!\!\!
  z^{-\frac{n}{2}}
  A_{n,r}(q)
  =
  \sum_{n\in 2\oZ+\overline{r+1}}\!\!
  z^{\frac{n}{2}}
  A_{n,r}(q)
  =
  \!\sum_{\substack{n=0\\ \nbar=\overline{r+1}}}^{2 k-1}\!\sum_{\ell\in\oZ}
  z^{k\ell + \frac{n}{2}}
  A_{2k\ell + n,r}.
\end{equation*}

Finally, the formula for $A_{n+2k\ell,r}(q)$ in~\bref{app:theF} allows
obtaining
\begin{multline*}
  \charW{+}{r}(q,z)
  =
  \ffrac{1}{\eta(q)}
  \sum_{\substack{n=0\\ \nbar=\overline{r+1}}}^{2 k-1}
  \Bigl(
  \hA_{n,r}(q)\theta_{n,k}(q,z)
  + \ffrac{2}{k}\,\hB_{n,r}(q)\theta'_{n,k}(q,z)
  + \ffrac{1}{k^2}\,
  \hC_{n,r}(q)\theta''_{n,k}(q,z)
  \!\Bigr)
  \\
  {}-
  \mfrac{z^{\frac{\rbar-1}{2}}}{\eta(q)^2}\,
  \sum_{\ell\in\oZ}\,z^{k(\ell+1)}
  q^{k(\ell+1)^2 + \ell(\overline{r+1})}
  \Bigl(\!
  \mfrac{2\ell\,q^{-k}}{
    1-z\,q^{2\ell}}\,\psi^+_{r}(q)
  -\mfrac{(2 \ell+1)q^{-\frac{k}{4}+\half(\overline{r+1})}}{
    1-z\,q^{2\ell+1}}\,\psi^-_{r}(q)
  \!\!\Bigr)
\end{multline*}
(again, with $\hC_{n,r}(q)$, $\hB_{n,r}(q)$, and $\hA_{n,r}(q)$
expressed as in~\eqref{ABC-help}), which readily yields the formula
in~\bref{chi-lemma}.

\begin{rem}\label{minus-reduced}
  It is easy to see that there is an equivalent representation for
  $\charW{-}{r}(q,z)$ and $\omega^-_r(q,z)$, with the Appell-function
  characteristics ``normalized'' to $\{0,1\}$:
  \begin{multline*}
    \charW{-}{r}(q,z)
    =-\ffrac{1}{\eta(q)}
    \!\!\!\!\sum_{\substack{n=0\\ \nbar=\overline{r+k+1}}}^{2 k-1}\!\!\!\!
    \Bigl(\!\!
    \bigl(\hA_{n-k,r}(q)-\ffrac{1}{4}\,\hC_{n-k,r}(q)\bigr)
    \theta_{n,k}(q,z)
    \\*[-8pt]
    \shoveright{{}
      +\ffrac{2}{k}\,\hB_{n-k,r}(q)\theta'_{n,k}(q,z)
      +\ffrac{1}{k^2}\,\hC_{n-k,r}(q)\theta''_{n,k}(q,z)
      \!\Bigr)}
    \\
    {}+\ffrac{\psi^+_{r}(q)}{\eta(q)^2}\;
    q^{-\frac{k}{4}}\Bigl(\ffrac{2}{k}\,\KKp_{\overline{k+r+1},k}(q,z,q^{-1})
    - \KK_{\overline{k+r+1},k}(q,z,q^{-1})
    \!\!\Bigr)
    -\ffrac{\psi^-_{r}(q)}{\eta(q)^2}
    \,\ffrac{2}{k}\,\KKp_{\overline{k+r+1},k}(q,z,1)
  \end{multline*}
  and
  \begin{multline*}
    \omega^-_r(q,z)
    =
    -\ffrac{1}{\eta(q)}
    \smash[t]{
      \!\!\!\!\sum_{\substack{n=0\\ \nbar=\overline{k+r+1}}}^{2 k-1}\!\!\!\!
      \Bigl(
      \hB_{n-k,r}(q)\theta_{n,k}(q,z)
      + \ffrac{1}{k}\,
      \hC_{n-k,r}(q)\theta'_{n,k}(q,z)
      \!\Bigr)}
    \\*[-4pt]
    {}+\ffrac{\psi^+_{r}(q)}{\eta(q)^2}\,
    q^{-\frac{k}{4}}\,\KK_{\overline{k+r+1},k}(q,z,q^{-1})
    -
    \ffrac{\psi^-_{r}(q)}{\eta(q)^2}
    \,\KK_{\overline{k+r+1},k}(q,z,1).
  \end{multline*}
\end{rem}

\section{Modular transformations}\label{sec:modular}
In this section, we use the decompositions in~\bref{chi-lemma}
and~\bref{lemma:minor} to derive modular transformation properties of
the functions~\eqref{the-set} occurring there as coefficients, among
which we are interested in the string functions $\hB_{n,r}(\tau)$ and
$\hA_{n,r}(\tau)$; that is, we prove the $S$-trans\-formation formulas
in~\bref{thm:B} and~\bref{thm:A}.

\subsection{$\smash{\boldsymbol{\hC_{n,r}(-\frac{1}{\tau})}}$}
For uniformity, we first rederive the well-known $S$-transformation of
the string functions $\hC_{n,r}(\tau)$.  {}From~\eqref{int-decomp},
\eqref{eq:theta-S}, and~\eqref{modular-eta},
\begin{align*}
  \charSL{}{r}(-\ffrac{1}{\tau},\ffrac{\nu}{\tau})
  &=-\ffrac{1}{\sqrt{2 k}}\,
  \ffrac{e^{i\pi k\frac{\nu^2}{2\tau}}}{
    \eta(\tau)}
  \sum_{\substack{m=0\\ \mbar=\overline{r+1}}}^{2 k-1}
  \hC_{m,r}(-\ffrac{1}{\tau})
  \sum_{n=0}^{2k-1}\!
  e^{-i\pi\frac{m\,n}{k}}\theta_{n, k}(\tau,\nu),
  \\
  \intertext{but in view of~\eqref{min-S-action} this is
    simultaneously equal to}
  &=-  \sqrt{\fffrac{2}{p}}\;
  \ffrac{e^{i\pi k\frac{\nu^2}{2\tau}}}{\eta(\tau)}\;  
  \sum_{s=1}^{p-1}\sin\!\ffrac{\pi r s}{p}\,
  \!\!\sum_{\substack{n=0\\ \nbar=\overline{s+1}}}^{2 k-1}\!\!
  \hC_{n,s}(\tau)
  \theta_{n,k}(\tau,\nu).
\end{align*}
Comparing the two expressions immediately yields~\eqref{modS-C}.

\subsection{$\smash{\boldsymbol{\hB_{n,r}(-\frac{1}{\tau})}}$}
\label{sec:S-B}
Following the same simple strategy to
find~$\hB_{n,r}(-\frac{1}{\tau})$ is somewhat more involved.  It is
technically convenient to introduce the linear combinations
\begin{align*}
  \Omega^a_r(\tau,\nu)
  &=\omega^+_r(\tau,\nu) + (-1)^a\omega^-_r(\tau,\nu),\\
  \KKK^a_\alpha(\tau,\nu,\mu)
  &=\KK_{\alpha,k}(\tau,\nu,\mu)+
  (-1)^ae^{-i\pi k\frac{\tau}{2} - i\pi k\mu}\KK_{\alpha+k,k}(\tau,\nu,\mu+\tau).
\end{align*}

\subsubsection{}\label{B-1st}From~\bref{lemma:minor},
\eqref{eq:theta-S}, \eqref{eq:theta'-S}, and~\bref{lemma:basic-S}, we
calculate
\begin{multline*}
  \Omega^a_r(-\ffrac{1}{\tau},\ffrac{\nu}{\tau})
  =\ffrac{2\,e^{i\pi k\frac{\nu^2}{2\tau}}}{\sqrt{2k}\,\eta(\tau)}
  \sum_{\substack{m=0\\ \mbar=\overline{r+1}}}^{2 k-1}
  \sum_{\substack{n=0\\ \nbar=\overline{a+1}}}^{2k-1}
  e^{-i\pi\frac{m\,n}{k}}
  \Bigl(
  \hB_{m,r}(-\ffrac{1}{\tau})
  \theta_{n,k}(\tau,\nu)
  \\[-6pt]
  \shoveright{{}+\ffrac{1}{k}\,\hC_{m,r}(-\ffrac{1}{\tau})
    \bigl(
    \tau\theta'_{n,k}(\tau,\nu) + \ffrac{k\nu}{2}\theta_{n,k}(\tau,\nu)
    \bigr)\!
    \Bigr)\ }
  \\[-3pt]
  %%% 
  {}+
  \ffrac{e^{i\pi k\frac{\nu^2}{2\tau}}\psi^+_{r}(-\frac{1}{\tau})}{
    i\eta(\tau)^2}\,
  \biggl[
  \KK_{\overline{a+1},k}(\tau,\nu,0) +
  (-1)^{r+1}e^{-i\pi k\frac{\tau}{2}}
  \KK_{\overline{a+1},k}(\tau,\nu,-\tau)
  \\[-3pt]
  \shoveright{{}+ e^{2i\pi \frac{k}{\tau}}
    \!\!\!\!\sum_{\beta\in\{0,1\}}
    \sum_{\substack{n=0\\ \nbar=\overline{a+1}}}^{2k-1}
    (-1)^{(r+1)\beta}
    \Phi(2k\tau,2k-n\tau-\beta k\tau)
    \theta_{n,k}(\tau,\nu)
    \biggr]}
  \\[-3pt]
  -\ffrac{e^{i\pi k\frac{\nu^2}{2\tau}}\psi^-_{r}(-\frac{1}{\tau})}{
    i\eta(\tau)^2}\,
  \biggl[
  (-1)^{a+1}\KK_{\overline{a+1},k}(\tau,\nu,0) +
  (-1)^{r+a+k}e^{-i\pi k\frac{\tau}{2}}
  \KK_{\overline{a+1},k}(\tau,\nu,-\tau)
  \\
  + e^{i\pi\frac{k}{2\tau}}
  \!\!\!\!\sum_{\beta\in\{0,1\}}\sum_{\substack{n=0\\ \nbar=\overline{a+1}}}^{2k-1}
  (-1)^{(r+1)\beta}
  \Phi(2k\tau, k-n\tau-\beta k\tau)
  \theta_{n,k}(\tau,\nu)
  \biggr].
\end{multline*}

\subsubsection{}\label{B-2nd}On the other hand, it follows
from~\bref{minor-modular} that
\begin{multline*}
  \Omega^a_r(-\ffrac{1}{\tau},\ffrac{\nu}{\tau})
  =
  \sqrt{\fffrac{2}{p}}\;e^{i\pi k\frac{\nu^2}{2\tau}}
  \biggl(
  \ffrac{i r}{2p}
  \Bigl[(-1)^r(1+(-1)^{a+p})\omega^+_p(\tau,\nu)
  +(1+(-1)^{a})\omega^-_p(\tau,\nu)\Bigr]
  \\
  \qquad\qquad\qquad\qquad\qquad{}
  +2\sum_{\substack{s=1\\ \sbar=a}}^{p-1}
  \Srstau
  \omega^+_s(\tau,\nu)
  + 2(-1)^{r}\!\!\sum_{\substack{s=1\\ \sbar=\overline{a+k}}}^{p-1}\!\!
  \Srstau
  \omega^-_s(\tau,\nu)
  \!\!\biggr)\\[-6pt]
  -\ffrac{2\nu}{\sqrt{2p}}\;e^{i\pi k\frac{\nu^2}{2\tau}}
  \sum_{\substack{s=1\\ \sbar=a}}^{p-1}
  \sin\!\ffrac{\pi r s}{p}\,\charSL{}{s}(\tau,\nu)
\end{multline*}
(see~\eqref{Srstau} for $\Srstau)$).  We next use the decompositions
of $\omega^{\pm}_s$ (and $\charSL{}{s}(\tau,\nu)$) again.  More
precisely, we express $\omega^+_s$ from~\bref{lemma:minor} and
$\omega^-_s$ from~\bref{minus-reduced}, which gives
\begin{small}%%%%%% KEEP THIS small!!
\begin{multline*}
  \Omega^a_r(-\ffrac{1}{\tau},\ffrac{\nu}{\tau})
  =
  \sqrt{\fffrac{2}{p}}\;
  \ffrac{e^{i\pi k\frac{\nu^2}{2\tau}}}{\eta(\tau)}
  \ffrac{i r}{2p}\times{}\\
  \shoveleft{{}\times\biggl[(-1)^r(1+(-1)^{a+k})
    \!\!\!\sum_{\substack{n=0\\ \nbar=\overline{a+1}}}^{2 k-1}\!\!\!
    \hB_{n,p}(\tau)\theta_{n,k}(\tau,\nu)
    -(1+(-1)^{a})
    \!\!\!\sum_{\substack{n=0\\ \nbar=\overline{a+1}}}^{2 k-1}\!\!\!
    \hB_{n-k,p}(\tau)\theta_{n,k}(\tau,\nu)}
  \\
  {}- (-1)^r(1+(-1)^{a+k})\Bigl(\!
  \ffrac{\psi^+_{p}(\tau)e^{-2i\pi k\tau}}{\eta(\tau)}\,
  \KK_{\overline{a+1},k}(\tau,\nu,{-2\tau})
  -
  \ffrac{\psi^-_{p}(\tau)e^{-i\pi\frac{k}{2}\tau}}{\eta(\tau)}\,
  \KK_{\overline{a+1},k}(\tau,\nu,{-\tau})
  \!\Bigr)
  \\  
  \shoveright{{}+(1+(-1)^{a})
    \Bigl(
    \ffrac{\psi^+_{p}(\tau)e^{-i\pi\frac{k}{2}\tau}}{\eta(\tau)}\,
    \KK_{\overline{a+1},k}(\tau,\nu,-\tau)
    -
    \ffrac{\psi^-_{p}(\tau)}{\eta(\tau)}\,\KK_{\overline{a+1},k}(\tau,\nu,0)
    \!\Bigr)\biggr]}
  \\
  \shoveleft{{}+2\sqrt{\fffrac{2}{p}}\;
    \ffrac{e^{i\pi k\frac{\nu^2}{2\tau}}}{\eta(\tau)}
    \sum_{\substack{s=1\\ \sbar=a}}^{p-1}
    \Srstau
    \biggl[\sum_{\substack{n=0\\ \nbar=\overline{a+1}}}^{2 k-1}\!\!
    \Bigl(
    \hB_{n,s}(\tau)\theta_{n,k}(\tau,\nu)
    +\ffrac{1}{k}\,\hC_{n,s}(\tau)\theta'_{n,k}(\tau,\nu)
    \!\Bigr)}
  \\
  \shoveright{
    -\ffrac{\psi^+_{s}(\tau)e^{-2i\pi k\tau}}{\eta(\tau)}\,
    \KK_{\overline{a+1},k}(\tau,\nu,-2\tau)
    +
    \ffrac{\psi^-_{s}(\tau)e^{-i\pi\frac{k}{2}\tau}}{\eta(\tau)}\,
    \KK_{\overline{a+1},k}(\tau,\nu,-\tau)
    \biggr]}
  \\
  \shoveleft{{}+ 2(-1)^{r}\sqrt{\fffrac{2}{p}}\;
    \ffrac{e^{i\pi k\frac{\nu^2}{2\tau}}}{\eta(\tau)}
    \sum_{\substack{s=1\\
        \sbar=\overline{a+k}}}^{p-1}\!\!
    \Srstau
    \biggl[-\!\!\sum_{\substack{n=0\\ \nbar=\overline{a+1}}}^{2 k-1}\!\!\!
    \Bigl(
    \hB_{n-k,s}(q)\theta_{n,k}(\tau,\nu)
    + \ffrac{1}{k}\,
    \hC_{n-k,s}(q)\theta'_{n,k}(\tau,\nu)
    \!\Bigr)}
  \\
  \shoveright{{}+
    \ffrac{\psi^+_{s}(\tau)e^{-i\pi\frac{k}{2}\tau}}{\eta(\tau)}\,
    \KK_{\overline{a+1},k}(\tau,\nu,-\tau)
    -
    \ffrac{\psi^-_{s}(\tau)}{\eta(\tau)}\,\KK_{\overline{a+1},k}(\tau,\nu,0)
    \biggr]}
  \\
  {}+
  \ffrac{2\nu}{\sqrt{2p}}\,\ffrac{e^{i\pi k\frac{\nu^2}{2\tau}}}{\eta(\tau)}
  \sum_{\substack{s=1\\ \sbar=a}}^{p-1}
  \sum_{\substack{n=0\\ \nbar=\overline{a+1}}}^{2 k-1}\!\!
  \sin\!\ffrac{\pi r s}{p}\,
  \hC_{n,s}(\tau)\theta_{n,k}(\tau,\nu).
\end{multline*}%
\end{small}%

\subsubsection{}We now compare the two expressions for
$\Omega^a_r(-\frac{1}{\tau},\frac{\nu}{\tau})$, in~\bref{B-1st}
and~\bref{B-2nd}.  The terms that explicitly involve $\nu$ already
coincide in view of~\eqref{modS-C}.  The terms involving $\theta'$ are
readily seen to coincide for the same reason (and because
of~\eqref{C-mirror}).

Next, comparing the terms involving $\KK\!\!\!$ (or, equivalently, the
residues of the two expressions for
$\Omega^a_r(-\frac{1}{\tau},\frac{\nu}{\tau})$), we recover the
transformations of the $(p,1)$-model characters $\psi^{\pm}_r(\tau)$
in~\eqref{eq:S-psiplus}--\eqref{eq:S-psiminus} (this seems to be a
remarkably complicated way to derive these simple formulas).  But most
importantly, some of the $\KK_{\overline{a+1},k}$-terms contribute to
$\theta_{n,k}$-terms in accordance with~\bref{sec:openqp}.  Comparing
the $\theta_{n,k}$-terms then gives the relation
\begin{equation}\label{BB-first}
  2\,\sqrt{2k}\,
  \hB_{m,r}(-\ffrac{1}{\tau})=
  \sum_{n=0}^{2k-1}e^{i\pi\frac{m\,n}{k}}
  \mathbb{B}_{r,b,n}(\tau),
  \qquad
  \begin{array}{l}
    \mbar=\overline{r+1},\\
    0\leq m\leq 2k\!-\!1,
  \end{array}
\end{equation}
where we temporarily use the notation
\begin{multline*}
  \mathbb{B}_{r,b,n}(\tau)
  =-\sqrt{\fffrac{2}{p}}
  \ffrac{i r}{2p}
  \Bigl[
  (-1)^r(1+(-1)^{n+1+k})
  \hB_{-n,p}(\tau)
  +(1+(-1)^{n+1})
  \hB_{n-k,p}(\tau)  
  \Bigr]
  \\
  {}-2\sqrt{\fffrac{2}{p}}
  \!\!\!\sum_{\substack{s=1\\ \sbar=\overline{n+1}}}^{p-1}\!\!
  \Srstau
  \hB_{-n,s}(\tau)
  - 2\sqrt{\fffrac{2}{p}}\;
  (-1)^{r}
  \!\!\!\!\!\sum_{\substack{s=1\\ \sbar=\overline{n+1+k}}}^{p-1}\!\!\!\!\!
  \Srstau
  \hB_{n-k,s}(\tau)
  \\
  {}+\ffrac{i\psi^+_{r}(-\frac{1}{\tau})e^{2i\pi \frac{k}{\tau}}}{\eta(\tau)}\,
  \!\!\!\!\sum_{\beta\in\{0,1\}}\!\!
  (-1)^{(r+1)\beta}
  \Phi(2k\tau,2k-n\tau-\beta k\tau)
  \\
  {}-\ffrac{i\psi^-_{r}(-\frac{1}{\tau})e^{i\pi\frac{k}{2\tau}}}{\eta(\tau)}\,
  \!\!\!\!\sum_{\beta\in\{0,1\}}\!\!
  (-1)^{(r+1)\beta}
  \Phi(2k\tau, k-n\tau-\beta k\tau).
\end{multline*}
We also used~\eqref{B-minus} here.

It now follows from~\eqref{B-open}, \eqref{B-minus},
and~\eqref{phi-tau} that
$\mathbb{B}_{r,b,n+k}(\tau)=(-1)^{r+1}\mathbb{B}_{r,b,n}(\tau)$, and
therefore Eq.~\eqref{BB-first} can be rewritten as\footnote{We
  simultaneously see that $0=\sum_{n=0}^{2k-1} e^{i\pi\frac{m\,n}{k}}
  \mathbb{B}_{r,b,n}(\tau)$ for $\mbar=\rbar$, which also follows from
  comparison of the $\theta_{n,k}$-terms above.}
\begin{equation*}
  \sqrt{2k}\,
  \hB_{m,r}(-\ffrac{1}{\tau})=
  \sum_{n=0}^{k-1}e^{i\pi\frac{m\,n}{k}}
  \mathbb{B}_{r,b,n}(\tau),
  \qquad
  \begin{array}{l}
    \mbar=\overline{r+1},\\
    0\leq m\leq 2k\!-\!1.
  \end{array}
\end{equation*}
But in the $\psi^{\pm}_r(-\frac{1}{\tau})$-terms in this sum, we then
have
\begin{multline*}
  \sum_{n=0}^{k-1}e^{i\pi\frac{m\,n}{k}}
  \ffrac{i\psi^+_{r}(-\frac{1}{\tau})}{\eta(\tau)}\,
  e^{2i\pi \frac{k}{\tau}}
  \!\!\!\!\sum_{\beta\in\{0,1\}}\!\!
  (-1)^{(r+1)\beta}
  \Phi(2k\tau,2k-n\tau-\beta k\tau)
  =\\[-4pt]
  =\sum_{n=0}^{2k-1}e^{i\pi\frac{m\,n}{k}}\,
  \ffrac{i\psi^+_{r}(-\frac{1}{\tau})}{\eta(\tau)}\,
  e^{2i\pi \frac{k}{\tau}}
  \Phi(2k\tau,2k-n\tau),
  \quad \mbar=\overline{r+1},\\
  \shoveleft{\mbox{}\kern-10pt\text{and subsequently
      using~\eqref{phi-2k} and then~\eqref{Phi-S}, we continue this
      as}}\\
  =\ffrac{i\psi^+_{r}(-\frac{1}{\tau})}{\eta(\tau)}\;
  e^{i\pi \frac{(m+2k)^2}{2 k\tau}}
  \Phi\bigl(\ffrac{\tau}{2k},1+\ffrac{m}{k}\bigr)
  =\sqrt{2 k}
  \ffrac{\psi^+_{r}(-\frac{1}{\tau})}{\sqrt{-i\tau}\eta(\tau)}\,
  \Phi\bigl(-\ffrac{2k}{\tau},-\ffrac{m}{\tau}\bigr).
\end{multline*}
Thus rewritten, this term (and the $\psi^{-}_r(-\frac{1}{\tau})$-term
similarly) naturally combines with the left-hand side
of~\eqref{BB-first} into
\begin{equation*}
  \sfBfirst_{m,r}(\tau)=\hB_{m,r}(\tau)
  -\ffrac{\psi^+_r(\tau)}{\eta(\tau)}\,\Phi(2k\tau, m\tau)
  +\ffrac{\psi^-_r(\tau)}{\eta(\tau)}\,\Phi(2k\tau, (m-k)\tau)
\end{equation*}
such that
\begin{multline*}
  \sfBfirst_{m,r}(-\ffrac{1}{\tau})=
  \\
  {}=
  -\ffrac{1}{\sqrt{k p}}\sum_{n=0}^{k-1}e^{i\pi\frac{m\,n}{k}}
  \biggl(\!\ffrac{i r}{2 p}
  \Bigl[(-1)^r(1+(-1)^{n+1+k})
  \hB_{-n,p}(\tau)
  +(1+(-1)^{n+1})
  \hB_{n-k,p}(\tau)  
  \Bigr]  
  \\[-4pt]
  {}+2
  \!\!\!\sum_{\substack{s=1\\ \sbar=\overline{n+1}}}^{p-1}\!\!
  \Srstau
  \hB_{-n,s}(\tau)
  +2 (-1)^{r}
  \!\!\!\!\!\sum_{\substack{s=1\\ \sbar=\overline{n+1+k}}}^{p-1}\!\!\!\!\!
  \Srstau
  \hB_{n-k,s}(\tau)\!\!\biggr)
\end{multline*}
for $\mbar=\overline{r+1}$ and $0\leq m\leq 2k\!-\!1$.
With~\eqref{B-eval} and after simple transformations, this can be
conveniently rewritten as
\begin{multline*}
  \sfBfirst_{m,r}(-\ffrac{1}{\tau})=
  \ffrac{(-1)^r}{\sqrt{k p}}
  \ffrac{i r}{2 p}
  \bigl(1+(-1)^{k+1}\bigr)
  \ffrac{\psi^+_p(\tau)}{2\eta(\tau)}  
  \\[-4pt]
  {}+\ffrac{1}{\sqrt{k p}}
  \sum_{\substack{s=1\\ \sbar=1}}^{p-1}
  \Srstau
  \ffrac{\psi^+_s(\tau)}{\eta(\tau)}
  + \ffrac{(-1)^{r}}{\sqrt{k p}}
  \!\!\sum_{\substack{s=1\\ \sbar=\overline{k+1}}}^{p-1}\!\!\!
  \Srstau
  \ffrac{\psi^-_s(\tau)}{\eta(\tau)}\qquad\qquad
  \\[-4pt]
  {}+\ffrac{(-1)^r}{\sqrt{k p}}
  \ffrac{2 r}{p}
  \!\!\!\sum_{\substack{n=1\\ \nbar=\overline{k+1}}}^{k-1}\!\!\!\!
  \sin\!\ffrac{\pi m\,n}{k}\,
  \hB_{-n,p}(\tau)
  -\ffrac{4 i}{\sqrt{k p}}
  \mathop{\sum_{n=1}^{k-1}
  \sum_{s=1}^{p-1}}\limits_{\overline{n+s}=1}
  \sin\!\ffrac{\pi m\,n}{k}\,
  \Srstau
  \hB_{-n,s}(\tau),
\end{multline*}
whence~\bref{thm:B} is immediate.

\subsection{$\smash{\boldsymbol{\hA_{n,r}(-\frac{1}{\tau})}}$}
\label{sec:S-A}
A similar calculation of $\hA_{n,r}(-\frac{1}{\tau})$ is
straightforward in principle but rather bulky in practical terms.  We
begin with introducing the linear combinations of characters
\begin{equation}
  \Xi^a_r(\tau,\nu)
  =\charW{+}{r}(\tau,\nu)
  + (-1)^a\bigl(\charW{-}{r}(\tau,\nu)
  + \ffrac{1}{4}\,\charSL{}{p-r}(\tau,\nu)\bigr)
\end{equation}

\subsubsection{}\label{A-1st}From~\bref{chi-lemma} and
\eqref{eq:theta-S}--\eqref{eq:theta''-S}, we
calculate
\begin{multline*}
  \Xi^a_r(-\ffrac{1}{\tau},\ffrac{\nu}{\tau})
  =\nu\Omega^a_r(-\ffrac{1}{\tau},\ffrac{\nu}{\tau})
  + \ffrac{2 e^{i\pi\frac{k \nu^2}{2\tau}}}{\sqrt{2 k}\eta(\tau)}
  \!\!\sum_{\substack{m=0\\ \mbar=\overline{r+1}}}^{2 k-1}
  \sum_{\substack{\ n=0\\ \nbar=\overline{a+1}}}^{2 k-1}\!\!\!
  e^{-i\pi\frac{m\,n}{k}}
  \biggl(
  \hA_{m,r}(-\ffrac{1}{\tau})\theta_{n, k}(\tau,\nu)
  \\
  +\ffrac{2\tau\!}{k}\,\hB_{m,r}(-\ffrac{1}{\tau})\theta'_{n, k}(\tau,\nu)
  + \ffrac{\tau^2\!}{k^2\!}\hC_{m,r}(-\ffrac{1}{\tau})
  \theta''_{n, k}(\tau,\nu)
  +\bigl(\ffrac{\tau}{4i\pi k}
  -\ffrac{\nu^2}{4}\bigr)
  \hC_{m,r}(-\ffrac{1}{\tau})
  \theta_{n, k}(\tau,\nu)
  \!\!\biggr)
  \\
  {}+\ffrac{\psi^+_{r}(-\frac{1}{\tau})}{i\tau\eta(\tau)^2}\,
  e^{2i\pi \frac{k}{\tau}}\Bigl(\ffrac{2}{k}\,
  {\KKK^{\overline{a+1}}_{\overline{r+1}}}'(
  -\ffrac{1}{\tau},\ffrac{\nu}{\tau},\ffrac{2}{\tau})
  - (2
  + \nu) %%%%%
  \KKK^{\overline{a+1}}_{\overline{r+1}}(
  -\ffrac{1}{\tau},\ffrac{\nu}{\tau},\ffrac{2}{\tau})\Bigr)
  \\
  {}- \ffrac{\psi^-_{r}(-\frac{1}{\tau})}{i\tau\eta(\tau)^2}\,
  e^{i\pi\frac{k}{2\tau}}\Bigl(\ffrac{2}{k}\,
  {\KKK^{\overline{a+1}}_{\overline{r+1}}}'(
  -\ffrac{1}{\tau},\ffrac{\nu}{\tau},\ffrac{1}{\tau})
  - (1
  + \nu) %%%%%
  \KKK^{\overline{a+1}}_{\overline{r+1}}(
  -\ffrac{1}{\tau},\ffrac{\nu}{\tau},\ffrac{1}{\tau})
  \!\Bigr)
\end{multline*}
We note that in the ``$\tau\nu\mu$'' notation,
\begin{equation*}
  \KKp_{\alpha,k}(\tau,\nu,\mu)=\ffrac{1}{2i\pi}
  \bigl(\ffrac{\dd}{\dd\nu} - \ffrac{\dd}{\dd\mu}\bigr)
  \KK_{\alpha,k}(\tau,\nu,\mu).
\end{equation*}
In substituting the $S$-transformed $\KKK^{a+1}_{\overline{r+1}}$
functions here, we evaluate the relevant combinations $e^{i\pi
  k\frac{\mu^2}{2\tau}}\bigl(\frac{2}{k}\,
{\KKK^{a+1}_{\overline{r+1}}}'(
-\frac{1}{\tau},\frac{\nu}{\tau},\frac{\mu}{\tau}) - (\nu +
\mu)\,\KKK^{a+1}_{\overline{r+1}}(
-\frac{1}{\tau},\frac{\nu}{\tau},\frac{\mu}{\tau})\bigr)$ at $\mu=1$
and~$2$ using the identity
\begin{multline*}
  \ffrac{1}{\tau}\,e^{i\pi k\frac{\mu^2}{2\tau}}
  \Bigl(
  \ffrac{2}{k}\,
  {\KKK^{a+1}_{\overline{r+1}}}'(
  -\ffrac{1}{\tau},\ffrac{\nu}{\tau},\ffrac{\mu}{\tau})
  - (\mu+\nu)\,\KKK^{a+1}_{\overline{r+1}}(
  \ffrac{1}{\tau},\ffrac{\nu}{\tau},\ffrac{\mu}{\tau})  
  \!\Bigr)
  =
  \\
  {}=\tau\,e^{i\pi k\frac{\nu^2}{2\tau}}\!
  \Bigl(\ffrac{2}{k}\,\KKp_{\overline{a+1},k}(\tau,\nu,\mu) +
  (-1)^{r+1}e^{-i\pi k\frac{\tau}{2} + i\pi k\mu}
  \ffrac{2}{k}\,\KKp_{\overline{a+1},k}(\tau,\nu,\mu\!-\!\tau)
  \\
  \shoveright{{}- (-1)^{r+1}e^{-i\pi k\frac{\tau}{2} + i\pi k\mu}
    \KK_{\overline{a+1},k}(\tau,\nu,\mu\!-\!\tau)\!\Bigr)}
  \\
  {}+ \ffrac{2}{k}\tau\,e^{i\pi k\frac{\nu^2 + \mu^2}{2\tau}}
  \!\!\!\!\sum_{\beta\in\{0,1\}}
  \sum_{\substack{n=0\\ \nbar=\overline{a+1}}}^{2k-1}\!\!
  (-1)^{(r+1)\beta}
  \Phi(2k\tau,k\mu\!-\!n\tau\!-\!\beta k\tau)
  \theta'_{n,k}(\tau,\nu)
  \\
  {}- \ffrac{\tau}{i\pi k}\,\ffrac{\dd}{\dd\mu}
  \Bigl(e^{i\pi k\frac{\nu^2 + \mu^2}{2\tau}}
  \!\!\!\!\sum_{\beta\in\{0,1\}}
  \sum_{\substack{n=0\\ \nbar=\overline{a+1}}}^{2k-1}\!\!
  (-1)^{(r+1)\beta}
  \Phi(2k\tau,k\mu\!-\!n\tau\!-\!\beta k\tau)
  \theta_{n,k}(\tau,\nu)\!\!\Bigr),
\end{multline*}
which readily follows from~\bref{lemma:basic-S}.  It is left to the
reader to substitute the last formula (twice) in the above expression
for $\Xi^a_r(-\fffrac{1}{\tau},\fffrac{\nu}{\tau})$.

\subsubsection{}\label{A-2nd}On the other hand, it follows
from~\bref{chi-modular} that
\begin{multline*}
  \Xi^a_{r}(-\ffrac{1}{\tau},\ffrac{\nu}{\tau})
  = \nu\,\Omega^a_r(-\ffrac{1}{\tau},\ffrac{\nu}{\tau}) +
  \\
  {}+\sqrt{\fffrac{2}{p}}\,
  e^{i\pi k\frac{\nu^2}{2\tau}}
  \biggl(\!
  \tau(1+(-1)^{k+a})\ffrac{(-1)^r i r}{2 p}\charW{+}{p}(\tau,\nu)
  {}+\tau(1+(-1)^{a})\ffrac{i r}{2 p}\charW{-}{p}(\tau,\nu)
  \\
  \begin{aligned}
    &+2\!\sum_{\substack{s=1\\ \sbar=\abar}}^{p-1}\!
    \tau\,\Srstau\charW{+}{s}(\tau,\nu)
    +2(-1)^r\!\!\sum_{\substack{s=1\\ \sbar=\overline{a+k}}}^{p-1}\!\!\!
    \tau\,\Srstau\charW{-}{s}(\tau,\nu)
    \\
    &+2\!\sum_{\substack{s=1\\ \sbar=\abar}}^{p-1}\!
    \ffrac{i r}{p}\,\tau\ffrac{s}{2p}\cos\!\ffrac{\pi r s}{p}\,
    \charSL{}{s}(\tau,\nu)
    + 
    2\sum_{\substack{s=1\\ \sbar=\abar}}^{p-1}
    \Bigl(
    \ffrac{\nu^2}{4}
    - \ffrac{s^2 \tau^2}{4p^2}\!\Bigr)
    \sin\!\ffrac{\pi r s}{p}\,
    \charSL{}{s}(\tau,\nu)
  \end{aligned}
  \\[-6pt]
  -\bigl(\ffrac{r^2}{2 p^2} + \ffrac{\tau}{2 i\pi p}\bigr)
  \sum_{\substack{s=1\\ \sbar=\abar}}^{p-1}
  \sin\!\ffrac{\pi r s}{p}\,\charSL{}{s}(\tau,\nu)
  \!\!\biggr).
\end{multline*}
We then use the decompositions for $\charW{\pm}{s}(\tau,\nu)$ (and
$\omega^{\pm}_s(\tau,\nu)$ and $\charSL{}{s}(\tau,\nu)$) again,
expressing $\charW{+}{s}$ from~\bref{chi-lemma} and $\charW{-}{s}$
from~\bref{minus-reduced}.  The substitution is totally
straightforward, but the result is rather cumbersome, and we leave it
to the reader to expand the last formula.

\subsubsection{}We next compare the two (rather cumbersome)
expressions for $\Xi^a_r(-\frac{1}{\tau},\frac{\nu}{\tau})$, resulting
from~\bref{A-1st} and~\bref{A-2nd}.  The terms proportional to $\nu$
are already written as $\nu\,\Omega^a_r$ and therefore cancel.  The
terms proportional to $\nu^2$ are readily seen to cancel due to the
$S$-transformation properties of~$\hC_{m,r}$.  The terms involving
$\theta''_{n,k}$ cancel for the same reason.  Further, all terms
involving $\theta'_{n,k}$ cancel due to the $S$-transformation
properties of $\hB_{m,r}$.\footnote{The calculation with
  $\charW{\pm}{r}$ alone establishes the transformations of
  $\hC_{m,r}$ and $\hB_{m,r}$ as well as of~$\hA_{m,r}$, but we prefer
  to have the formula for $\hB_{m,r}(-\frac{1}{\tau})$ already derived
  and to use it in the (rather tedious) $\charW{\pm}{r}$-calculation
  only for control.}  After cancellations of the
$\KKp_{\overline{a+1},k}$ and $\KK_{\overline{a+1},k}$, based on the
identity
\begin{multline*}
  e^{-2i\pi k\tau}
  \Bigl(\ffrac{2}{k}\;\KKp_{\overline{a+1},k}(\tau,\nu,-2\tau)
  - 2 \KK_{\overline{a+1},k}(\tau,\nu,-2\tau)\!\Bigr)=\\
  {}= \ffrac{2}{k}\;
  \KKp_{\overline{a+1},k}(\tau,\nu,0)  
  + \!\!\sum_{\substack{n=0\\ \nbar=\overline{a+1}}}^{2k-1}\!\!\!
  \Bigl(  
  \ffrac{n}{k}\;
  e^{-i\pi\tau\frac{n^2}{2k}}\theta_{n,k}(\tau,\nu)
  - \ffrac{2}{k}\,e^{-i\pi\tau\frac{n^2}{2k}}\theta'_{n,k}(\tau,\nu)
  \!\Bigr)
\end{multline*}
following from~\bref{sec:openqp}, we obtain
\begin{equation*}
  2\sqrt{2k}  
  \hA_{m,r}(-\ffrac{1}{\tau})
  = \sum_{n=0}^{2 k-1} e^{i\pi\frac{m\,n}{k}} \mathbb{A}_{r,b,n}(\tau),
  \quad
  \begin{array}{l}
    \mbar=\overline{r+1},\\
    0\leq m\leq 2k-1,
  \end{array}
\end{equation*}
where we introduce the temporary notation
\begin{multline*}
  \mathbb{A}_{r,b,n}(\tau)=\sqrt{\fffrac{2}{p}}
    \Biggl(
    \tau(1+(-1)^{k+n+1})\ffrac{(-1)^r i r}{2 p}
    \hA_{-n,p}(\tau)
    {}-\tau(1+(-1)^{n+1})\ffrac{i r}{2 p}
    \hA_{n-k,p}(\tau)
  \\
  +2\!\sum_{\substack{s=1\\ \sbar=\overline{n+1}}}^{p-1}\!
  \tau\,\Srstau
  \hA_{-n,s}(\tau)
  - 
  2(-1)^r
  \!\!\!\!\!\sum_{\substack{s=1\\ \sbar=\overline{n+k+1}}}^{p-1}\!\!\!\!
  \tau\,\Srstau
  \hA_{n-k,s}(\tau)
  \\
  +\!\sum_{\substack{s=1\\ \sbar=\overline{n+1}}}^{p-1}\!
  \ffrac{i r}{p}\,\tau\ffrac{p -2 s}{2 p}
  \cos\!\ffrac{\pi r s}{p}\,  
  \hC_{n,s}(\tau)
  + \sum_{\substack{s=1\\ \sbar=\overline{n+1}}}^{p-1}\!\!
  \Bigl(
  \ffrac{(s^2 - p s) \tau^2}{2 p^2}
  + \ffrac{r^2}{2 p^2}
  - \ffrac{\tau}{i\pi p k}
  \Bigr)
  \sin\!\ffrac{\pi r s}{p}\,
  \hC_{n,s}(\tau)
  \!\!\Biggr)
  \\
  {}+\ffrac{\psi^+_{r}(-\frac{1}{\tau})}{i\eta(\tau)}\,
  \ffrac{\tau}{i\pi k}\,\ffrac{\dd}{\dd\mu}
  \!\!\sum_{\beta\in\{0,1\}}
  (-1)^{(r+1)\beta}
  e^{i\pi k\frac{\mu^2}{2\tau}}
  \Phi(2k\tau,k\mu\!-\!n\tau\!-\!\beta k\tau)
  \Bigr|_{\mu=2}
  \\
  {}- \ffrac{\psi^-_{r}(-\frac{1}{\tau})}{i\eta(\tau)}\,
  \ffrac{\tau}{i\pi k}\,\ffrac{\dd}{\dd\mu}
  \!\!\sum_{\beta\in\{0,1\}}
  (-1)^{(r+1)\beta}
  e^{i\pi k\frac{\mu^2}{2\tau}}
  \Phi(2k\tau,k\mu\!-\!n\tau\!-\!\beta k\tau)
  \Bigr|_{\mu=1}.
\end{multline*}
We also used~\eqref{A-minus} here.  It now follows
from~\eqref{A-open}, \eqref{A-minus}, and~\eqref{phi-tau} that
$\mathbb{A}_{r,b,n+k}(\tau)=(-1)^{r+1}\mathbb{A}_{r,b,n}(\tau)$, and
therefore
\begin{equation*}
  \sum_{n=0}^{2 k-1} e^{i\pi\frac{m\,n}{k}} \mathbb{A}_{r,b,n}(\tau)=
  2 \sum_{n=0}^{k-1} e^{i\pi\frac{m\,n}{k}} \mathbb{A}_{r,b,n}(\tau).  
\end{equation*}
But in the $\psi^{\pm}_r(-\frac{1}{\tau})$-terms in the sum in the
right-hand side, we then have (see~\bref{sec:Phi}
and~\bref{sec:Phi-S})
\begin{multline*}
  \sum_{n=0}^{k-1} e^{i\pi\frac{m\,n}{k}}
  \ffrac{\psi^{\pm}_{r}(-\frac{1}{\tau})}{i\eta(\tau)}\,
  \ffrac{\tau}{i\pi k}\,\ffrac{\dd}{\dd\mu}
  \!\sum_{\beta\in\{0,1\}}\!
  (-1)^{(r+1)\beta}
  e^{i\pi k\frac{\mu^2}{2\tau}}\,
  \Phi(2k\tau,k\mu\!-\!n\tau\!-\!\beta k\tau)
  =
  \\
  =
  \ffrac{\tau}{i\pi k}\,\ffrac{\dd}{\dd\mu}
  \ffrac{\psi^{\pm}_{r}(-\frac{1}{\tau})}{i\eta(\tau)}\,
  e^{i\pi k\frac{\mu^2}{2\tau}}\,
  \sum_{n=0}^{2 k-1} e^{i\pi\frac{m\,n}{k}}
  \Phi(2k\tau,k\mu\!-\!n\tau)
  =
  \\
  =\ffrac{\tau}{i\pi k}\,\ffrac{\dd}{\dd\mu}
  \ffrac{\psi^{\pm}_{r}(-\frac{1}{\tau})}{i\eta(\tau)}\,
  e^{i\pi\frac{(m + k\mu)^2}{2 k\tau}}
  \Phi\bigl(\ffrac{\tau}{2k},\ffrac{\mu}{2} + \ffrac{m}{2 k}\bigr)
  =\\
  =\sqrt{2k}\,
  \ffrac{2\psi^{\pm}_{r}(-\frac{1}{\tau})}{\sqrt{-i\tau}\eta(\tau)}\,
  \Phi'\bigl(-\ffrac{2k}{\tau},-\ffrac{m+(\mu-2)k}{\tau}\bigr).
\end{multline*}
Hence, defining $\sfA_{m,r}(\tau)$ as in~\eqref{sfA-def}, we obtain
the $S$-transform formula in~\bref{thm:A}.

\section{Conclusions}
The higher string functions $\hA_{n,r}(\tau)$ and $\hB_{n,r}(\tau)$
are not ``arbitrary'' analogues of the $\hC_{n,r}(\tau)$: there is an
underlying representation-theory picture described in~\cite{[S-sl2]}.
The associated conformal field theory construction (the $W$-algebra
in~\cite{[S-sl2]}) may then be considered the rationale for the higher
string functions to have interesting modular properties.

A ``feedback'' of modular transformations to conformal field theory is
that they come to play the role of a strong consistency check (e.g.,
for the field content) whenever representation-theory details are not
known, as is the case with the logarithmic extension of the
parafermion theory, where only the characters are available but the
field-theory picture is presently obscure.  As in the previously known
cases of logarithmic $(p,1)$ and $(p,p')$
models~\cite{[FGST],[FGST3]}, the degree of the polynomials in $\tau$
occurring in modular transformations may be expected to correlate with
the Jordan cell sizes in indecomposable representations of the
corresponding extended algebra, but the representation-theory
interpretation of the occurrences of the ``$\Phi$-constants''
$\Phi(2k\tau, m\tau)$ (times the $(p,1)$-model characters) is a
challenging problem.
%% If a consistent logarithmic parafermion theory can be constructed
%% directly, it will be an interesting extension of the $(p,1)$
%% logarithmic models.

Modular transformations are related to fusion, via the Verlinde
formula in rational conformal field theories~\cite{[V]} and via its
generalizations in logarithmic theories~\cite{[FHST],[FK],[GR2]};
whether the modular transformations derived in this paper lead to any
reasonable nonsemisimple fusion algebra remains an interesting
problem.

The celebrated form of $\hC_{n,r}(q)$ first found in~\cite{[LP]} has
been the subject of considerable attention since then; it would be
interesting to find an extension of such representations to the
``logarithmically extended parafermionic characters'' $\hA_{m,r}(q)$
and $\hB_{m,r}(q)$.  Different ``fermionic-type'' character formulas
may also be mentioned in this connection (see \cite{[KKMcCM]} and the
numerous subsequent papers, in particular,
e.g.,~\cite{[BLS],[ANOT+],[ABD]} and the references therein).  Their
counterparts for $\hA_{m,r}(q)$ and $\hB_{m,r}(q)$ may also be
interesting.

\subsubsection*{Acknowledgments} I am grateful to A.~Gainutdinov for
the useful comments.  This paper was supported in part by the RFBR
grant 07-01-00523 and grant LSS-4401.2006.2.

\appendix
\section{Theta-function conventions}\label{app:theta}
The level-$\varkappa$ theta-functions are defined as
\begin{gather}\label{thetaKac}
  \theta_{r, \varkappa}(q,z)
  =\smash[b]{\sum_{\iota\in\oZ+\frac{r}{2 \varkappa}}}
  q^{\varkappa \iota^2}z^{\varkappa \iota}.
\end{gather}
We set
\begin{gather}
  \theta'_{r, \varkappa}(q,z)
  \smash{=z\ffrac{\dd}{\dd z}\,\theta_{r, \varkappa}(q,z),
    \quad
    \theta''_{r, \varkappa}(q,z)=\Bigl(z\ffrac{\dd}{\dd z}\Bigr)^2}
  \theta_{r, \varkappa}(q,z)\\[-4pt]
  \intertext{and}
  \theta'_{r, \varkappa}(q)=\theta'_{r, \varkappa}(q,z)\Bigr|_{z=1}.
\end{gather}

The quasiperiodicity properties of theta-functions are expressed as
\begin{align}\label{theta-sf}
  \theta_{r,\varkappa}(q,z q^n) &=
  q^{-\varkappa\frac{n^2}{4}}z^{-\varkappa\frac{n}{2}}\,
  \theta_{r+\varkappa n, \varkappa}(q,z),
\end{align}
with $\theta_{r + \varkappa n, \varkappa}(q,z)=\theta_{r,
  \varkappa}(q,z)$ for \textit{even} $n$.

The modular $T$-transform of theta-functions
is
\begin{align}
  \theta_{r, \varkappa}(\tau+1, \nu)
  &=e^{i\pi\frac{r^2}{2\varkappa}}\,\theta_{r, \varkappa}(\tau, \nu)\\
  \intertext{and the $S$-transform is}
  \label{eq:theta-S}
  \theta_{r, \varkappa}(-\ffrac{1}{\tau}, \ffrac{\nu}{\tau})
  &=e^{i\pi\frac{\varkappa\nu^2}{2\tau}}
  \sqrt{\fffrac{-i\tau}{2 \varkappa}}\,\sum_{s=0}^{2\varkappa-1}\!
  e^{-i\pi\frac{r s}{\varkappa}}\theta_{s, \varkappa}(\tau,\nu).
  \\[-4pt]
  \intertext{Therefore,}
  \label{eq:theta'-S}
  \theta'_{r, \varkappa}(-\ffrac{1}{\tau}, \ffrac{\nu}{\tau})
  &=e^{i\pi\frac{\varkappa\nu^2}{2\tau}}
  \smash[t]{\sqrt{\fffrac{-i\tau}{2 \varkappa}}
    \,\sum_{s=0}^{2\varkappa-1}\!  e^{-i\pi\frac{r s}{\varkappa}}}
  \Bigl(
  \tau\theta'_{s, \varkappa}(\tau,\nu)
  + \ffrac{\varkappa\nu}{2}\,\theta_{s, \varkappa}(\tau,\nu)
  \Bigr),
  \\
  \label{eq:theta''-S}
  \theta''_{r, \varkappa}(-\ffrac{1}{\tau}, \ffrac{\nu}{\tau})
  &=e^{i\pi\frac{\varkappa\nu^2}{2\tau}}
  \sqrt{\fffrac{-i\tau}{2 \varkappa}}\,\sum_{s=0}^{2\varkappa-1}\!
  e^{-i\pi\frac{r s}{\varkappa}}
  \Bigl(
  \tau^2\theta''_{s, \varkappa}(\tau,\nu)
  + \varkappa\nu\tau\theta'_{s, \varkappa}(\tau,\nu)\\*[-4pt]
  &\qquad\qquad\qquad\qquad\qquad\qquad\qquad\quad
  {}+\bigl(\ffrac{\varkappa^2\nu^2}{4}+\ffrac{\varkappa\tau}{4i\pi}\bigr)
  \theta_{s, \varkappa}(\tau,\nu)
  \!\Bigr).\notag
\end{align}
The price paid for abusing notation is that
$\theta'_{r,\varkappa}(\tau,\nu)=\fffrac{1}{2i\pi}\fffrac{\dd}{\dd\nu}
\theta_{r,\varkappa}(\tau,\nu)$.

We also use the Jacobi theta-functions
\begin{gather}\label{vartheta}
  \vartheta_{1,1}(q, z)
  =  \sum_{m\in\oZ}q^{\half(m^2 - m)} (-z)^{-m}
  = \prod_{m\geq0}\!(1\!-\!z^{-1} q^m)
  \prod_{m\geq1}\!(1\!-\!z q^m)\prod_{m\geq1}\!(1\!-\!q^m),\\
  \label{plain-theta}
  \vartheta(q,z)
  =\sum_{m\in\oZ}q^{\frac{m^2}{2}}z^m
\end{gather}
related to~\eqref{thetaKac} as
\begin{align*}
  \theta_{r, \varkappa}(q,z)
  &=z^{\frac{r}{2}}\,q^{\frac{r^2}{4 \varkappa}}\,
  \vartheta(q^{2 \varkappa},z^{\varkappa} q^r).
\end{align*}
For the function
\begin{equation}\label{Omega}
  \Omega(q,z)=q^{\frac{1}{8}}z^{\half}\vartheta_{1,1}(q,z),
\end{equation}
we then have the $S$-transformation formula
\begin{align*}
  \Omega(-\ffrac{1}{\tau},\ffrac{\nu}{\tau})
  &=-i\sqrt{-i\tau}\,e^{i\pi\frac{\nu^2}{2\tau}}\Omega(\tau,\nu).
\end{align*}

The eta function
\begin{gather}\label{eta}
  \eta(q)  
  =q^{\frac{1}{24}}
  \smash{\prod\limits_{m=1}^\infty(1-q^{m})}
\end{gather}
transforms as
\begin{gather}\label{modular-eta}
  \eta(\tau+1)=\smash{e^{\frac{i\pi}{12}}\eta(\tau),\qquad
    \eta(-\ffrac{1}{\tau})=
    \sqrt{-i\tau}\,\eta(\tau)}.
\end{gather}

\section{Properties of the Appell functions~\cite{[STT]}}
\label{sec:Appell}
\subsection{Definition and simple properties}
\subsubsection{} For $k\in\oN$, the level-$k$ Appell function is
defined as~\cite{[STT]}
\begin{gather}
  \K_k(q,x,y)=
  \sum_{m\in\oZ}\mfrac{\displaystyle
    q^{\frac{m^2 k}{2}} x^{m k}}{\displaystyle
    1 - x\,y\, q^{m}},
\end{gather}
It has a number of properties that nontrivially generalize the
theta-function properties, the crucial ones being the ``open
quasiperiodicity''
\begin{equation*}
  \K_k(q^2,xq^{-\frac{2n}{k}},yq^{\frac{2n}{k}})
  = (xy)^{n}\K_k(q^2,x, y) +
%%   \begin{cases}
%%     \displaystyle
    \sum_{r=1}^{n}(xy)^{n-r}\,x^r\,q^{-\frac{r^2}{k}}\theta_{-2r,k}(q,x),
%%     &
    \quad n\in\oN
%%     ,\\
%%     \displaystyle
%%     -\!\sum_{r=n+1}^{0}(x y)^{n-r}\,x^r\,q^{-\frac{r^2}{k}}\theta_{-2r,k}(q,x),
%%     \kern-4pt
%%     &n\in-\oN
%%   \end{cases}\kern-6pt
\end{equation*}
(where it is worth noting that the $x$ and $y$ variables separate in
the extra terms), and the modular transformation properties,
Eq.~\eqref{Knew-Stransf-fin} in particular (where, remarkably, $\nu$
and $\mu$ also separate in the extra terms in the right-hand side).

\subsubsection{}
In this paper, we need a version of $\K_k$ with the ``period'' $q^2$
and with characteristics.  We define
\begin{equation}\label{eq:KK}
    \KK_{\alpha,k}(q,x,y)={}
    q^{-\frac{\alpha^2}{4k}}\,y^{\frac{\alpha}{2}}
    \sum_{m\in\oZ+\frac{\alpha}{2k}}\frac{q^{k m^2}\,x^{k m}}{
      1 - x\,y\,q^{2 m - \frac{\alpha}{k}}}
%%   ={}&(x\,y)^{\alpha/2}\sum_{m\in\oZ}
%%   \frac{q^{k m^2 + \alpha m}\,x^{k m}}{1 - x\,y\,q^{2m}}\\
    {}={}
    (x\,y)^{\frac{\alpha}{2}}\,
    \K_k(q^2,x\,q^{\frac{\alpha}{k}},y\,q^{-\frac{\alpha}{k}}).
\end{equation}
The open quasiperiodicity relation above
implies that
\begin{equation}\label{K-char-reduce}
  \KK_{\alpha+2n,k}(q,x,y)=
  \KK_{\alpha,k}(q,x,y)
%%   -\sum_{b=0}^{n-1}q^{-\frac{1}{4k}(\alpha + 2b)^2}\,
%%   y^{\frac{\alpha}{2} + b}\,\theta_{\alpha+2b,k}(q,x),
    -\sum_{\substack{m=0\\ \mbar=\overline{\alpha}}}^{2 n-1}q^{-\frac{m^2}{4k}}\,
  y^{\frac{m}{2}}\,\theta_{m,k}(q,x),
\end{equation}
and therefore the characteristic $\alpha$ in~\eqref{eq:KK} can be
reduced modulo~$2$ at the expense of theta functions.

\subsubsection{Open quasiperiodicity in the third argument}\label{sec:openqp}
It follows from the formulas in~\cite{[STT]} (or can be easily derived
from the definition) that $\KK_{\alpha,k}$ satisfies an open
quasiperiodicity property with respect to the third argument:
\begin{align*}
  \KK_{\alpha,k}(q,x,y\,q^{-2n})&=
  q^{k n^2}\,y^{-k n}\Bigl(
  \KK_{\alpha,k}(q,x,y)
  - \sum_{b=0}^{n k -1}q^{-\frac{1}{4k}(\alpha+2 b)^2}
  \,y^{\frac{\alpha}{2} + b}\,\theta_{\alpha+2b,k}(q,x)
  \!\Bigr)
%%   ,
%%   \\
%%   \KK_{\alpha,k}(q,x,y\,q^{2n})&=
%%   q^{k n^2}\,y^{k n}\Bigl(
%%   \KK_{\alpha,k}(q,x,y)
%%   + \sum_{b=1}^{n k}q^{-\frac{1}{4k}(\alpha-2 b)^2}
%%   \,y^{\frac{\alpha}{2} - b}\,\theta_{\alpha-2b,k}(q,x)
%%   \!\Bigr)
\end{align*}
for $n\in\oN$.  \ In the text, we use this formula in the form
\begin{equation*}
  \KK_{\overline{a+1},k}(\tau,\nu,\mu\!-\!2\tau)=
  e^{2i\pi k\tau - 2i\pi k\mu}
  \Bigl(
  \KK_{\overline{a+1},k}(\tau,\nu,\mu)
  -\!\!\sum_{\substack{n=0\\ \nbar=\overline{a+1}}}^{2k-1}\!\!\!
  e^{-i\pi\tau\frac{n^2}{2k} + i\pi\mu n}\theta_{n,k}(\tau,\nu)
  \!\Bigr).
\end{equation*}
%% Hence, in particular,
%% \begin{multline*}
%%   \KKp_{\overline{a+1},k}(\tau,\nu,\mu-\!2\tau)
%%   - k \KK_{\overline{a+1},k}(\tau,\nu,\mu-\!2\tau)
%% %%   =\ffrac{1}{2i\pi}\bigl(\ffrac{\dd}{\dd\nu} - \ffrac{\dd}{\dd\mu}\bigr)
%% %%   \KK_{\overline{a+1},k}(\tau,\nu,\mu\!-\!2\tau)
%% %%   - k \KK_{\overline{a+1},k}(\tau,\nu,\mu-\!2\tau)\\
%%   = e^{2i\pi k\tau - 2i\pi k\mu}
%%   \KKp_{\overline{a+1},k}(\tau,\nu,\mu)  
%%   \\*
%%   {}+ e^{2i\pi k\tau - 2i\pi k\mu}
%%   \sum_{\substack{n=0\\ \nbar=\overline{a+1}}}^{2k-1}\!\!\!
%%   \Bigl(  
%%   \ffrac{n}{2}\;e^{-i\pi\tau\frac{n^2}{2k} + i\pi\mu n}\theta_{n,k}(\tau,\nu)
%%   - e^{-i\pi\tau\frac{n^2}{2k} + i\pi\mu n}\theta'_{n,k}(\tau,\nu)
%%   \!\Bigr).
%% \end{multline*}

\subsubsection{Period-changing formula}
We  recall the elementary theta-func\-tion identity
\begin{gather}\label{theta-rewrite}
  \vartheta(q,z)=
  \sum_{s=0}^{u-1}
  q^{\frac{s^2}{2}}\,z^s
  \vartheta(q^{u^2},z^u q^{us}),\quad u\in\oN.
\end{gather}
Similarly to~\eqref{theta-rewrite}, we relate $\K_{k}(q^2,x,y)$ to
suchlike functions with the ``period'' $q^{u^2}$ for any $u\in\oN$ as
\begin{equation*}
  \K_{k}(q^2,x,y)
  =\sum_{a=0}^{u-1}\sum_{b=0}^{u-1}q^{-k a^2}\,y^{-k a}\,
  \KK_{\frac{2 b}{u} + \frac{2 k a}{u},k}(q^{u^2}, x^{u}, y^{u}\,q^{2 u a}).
\end{equation*}
This formula may not be very useful for general $u$ because of the
fractional characteristics in the right-hand side, but for $u=2$ it
takes the simple form
\begin{equation}\label{u2}
  \K_{k}(q^{\half},x^{\half},y^{\half})
%%   =\!\sum_{a,b\in\{0,1\}}\!q^{-\frac{k a^2}{4}}\,y^{-\frac{k a}{2}}\,
%%   \KK_{b + k a, k}(q, x,y\,q^a).
  =\!\!\!\sum_{\gamma,\beta\in\{0,1\}}\!\!\!
  q^{-\frac{k \gamma^2}{4}}\,y^{-\frac{k \gamma}{2}}\,
  \KK_{\beta + k \gamma, k}(q, x,y\,q^{\gamma})
  =\!\sum_{\beta\in\{0,1\}}\!
  \KKK^0_\beta(q,x,y).
\end{equation}

\subsection{The $\Phi$ function}\label{sec:Phi}
The $\Phi$ function defined in~\eqref{Phi-def} can be equivalently
written as the $b$-cycle integral~\cite[Eq.~(A.5)]{[STT]}
\begin{equation}
  \sqrt{-i\tau}\,\Phi(\tau,\mu)
  = i \int_{0}^{\tau} d\lambda\,
  \EX{i\pi\frac{\lambda^2 - 2\lambda\mu}{\tau}}\,
  \K_1(\tau,\lambda+\varepsilon-\mu,\mu)%
\end{equation}%
(where an infinitesimal positive real $\varepsilon$ specifies the
prescription to bypass the singularities).  This close relative of the
Mordell integral has a number of useful
properties~\cite{[STT]}. First, it is ``open-double-quasiperiodic''
under shifts of the argument by $n+m\tau$, $m,n\in\oZ$:
\begin{align}\label{phi-1}
%%   &\begin{aligned}
    \Phi(\tau, \mu\!+\!n) &=
    \EX{-i\pi\frac{n^2}{\tau}
      \!-\!2i\pi n\frac{\mu}{\tau}}\,\Phi(\tau,\mu)
    + \ffrac{i}{\sqrt{-i\tau}}
    \sum_{j=1}^n
    \EX{i\pi\frac{j(j-2 n)}{\tau}
      \!-\!2i\pi j\frac{\mu}{\tau}},
%%     \\
%%     \Phi(\tau,\mu\!-\!n) &=
%%     \EX{-i\pi\frac{n^2}{\tau}
%%       \!+\!2i\pi n\frac{\mu}{\tau}}\,\Phi(\tau,\mu)
%%     - \ffrac{i}{\sqrt{-i\tau}}    
%%     \sum_{j=0}^{n-1}
%%     \EX{i\pi\frac{j(j-2 n)}{\tau}
%%       \!+\!2i\pi j\frac{\mu}{\tau}},
%%   \end{aligned}
  \quad n\in\oN,
  \\[-4pt]
%%   \intertext{and}
  \label{phi-tau}
%%   &\begin{aligned}
%%     \Phi(\tau,\mu\!+\!m\tau) &= \Phi(\tau,\mu)
%%     - \smash[t]{\sum_{j=1}^{m}
%%       \ex{-i\pi\frac{(\mu+j\tau)^2}{\tau}}},\\
    \Phi(\tau,\mu\!-\!m\tau) &= \Phi(\tau,\mu)
    + \sum_{j=0}^{m-1}
    \ex{-i\pi\frac{(\mu - j\tau)^2}{\tau}},
%%   \end{aligned}
  \quad m\in\oN.
\end{align}
Together with the ``reflection'' property
\begin{align}\label{Phi-reflect}
  \Phi(\tau,-\mu) &=
  \ffrac{-i}{\sqrt{-i\tau}}
  - \EX{-i\pi\frac{\mu^2}{\tau}} - \Phi(\tau,\mu),
%%   \\
%%   \intertext{or equivalently,}
%%   \label{Phi-reflect2}
%%   \Phi(\tau,\mu+1) &=
%%   \smash[t]{-\EX{-2i\pi\frac{\mu + \half}{\tau}}} \Phi(\tau,-\mu-\tau),
\end{align}
this allows evaluating $\Phi(\tau,\,{\cdot}\,)$ at some (not all) of
the half-period points:
\begin{alignat*}{2}
%%   \label{explicit-phi-dual}
%%   &\begin{aligned}
  \Phi(\tau,\ffrac{n}{2})&=
  -\half\;\ex{-i\pi\frac{n^2}{4\tau}}
  + \ffrac{i}{2\sqrt{-i\tau}}\sum_{j=1}^{n-1}
  \ex{i\pi\frac{j(j-n)}{\tau}},\quad& n&{}\geq1,
%%     \\
%%     \Phi(\tau,-\ffrac{n}{2})&=
%%     -\half\;\ex{-i\pi\frac{n^2}{4\tau}}
%%     -\ffrac{i}{2\sqrt{-i\tau}}\sum_{j=0}^{n}
%%     \ex{i\pi\frac{j(j-n)}{\tau}},
%%     \quad& n&{}\geq0,
%%   \end{aligned}
  \\[-6pt]
  \intertext{and}
%%   \label{explicit-phi}
%%   &\begin{aligned}
  \Phi(\tau,\ffrac{m\tau}{2})&=
  \smash[t]{-\ffrac{i}{2\sqrt{-i\tau}}
    -\half \sum_{j=0}^m
    \ex{-i\pi\tau\frac{(m-2j)^2}{4}}},\quad& m&{}\geq0,
%%   \\
%%   \Phi(\tau,-\ffrac{m\tau}{2})&=-\ffrac{i}{2\sqrt{-i\tau}}
%%   + \half \sum_{j=1}^{m-1}
%%   \ex{-i\pi\tau\frac{(m-2j)^2}{4}},
%%   \quad& m&{}\geq1,
%%   \end{aligned}
\end{alignat*}
whence $\Phi(\tau,\fffrac{n}{2} + \fffrac{m}{2}\,\tau)$ \textit{with
  even $n m$} follow via the open quasiperiodicity formulas above
(formulas~\eqref{phi-1}--\eqref{Phi-reflect} do \textit{not} allow
finding $\Phi(\tau,\fffrac{n}{2} + \fffrac{m}{2}\,\tau)$ with $n$ and
$m$ both odd).

Next, a simple ``scaling law'' is given by
\begin{gather}\label{phi-u}
  \Phi(\tau,\mu) = \sum_{b=0}^{u-1}\Phi(u^2 \tau, u \mu - b u
  \tau),\qquad u\in\oN.
\end{gather}
In the case of ``scaling'' with an even factor, we have the
formula~\cite{[STT]}
\begin{gather}\label{phi-2k}
  \sum_{n=0}^{2k-1}\EX{i\pi\frac{m\,n}{k}}
  \Phi(2k\tau,2k\mu\!-\!n\tau)
  =\EX{i\pi\frac{[m]_{2k}^2}{2k\tau}
    \!+\!2i\pi\frac{\mu}{\tau}[m]_{2k}}
  \Phi(\ffrac{\tau}{2k},\mu\!+\!\ffrac{[m]_{2k}}{2k}).
%%   \\
%%   \ffrac{1}{2k}\sum_{n=0}^{2k-1}
%%   \EX{
%%     2i\pi\frac{n\mu}{\tau}
%%     \!-\!i\pi\frac{m\,n}{k}
%%     \!+\!i\pi\frac{n^2}{2k\tau}}
%%   \Phi(\ffrac{\tau}{2k},\mu\!+\!\ffrac{n}{2k})
%%   =\Phi(2k\tau,2k\mu\!-\![m]_{2k}\,\tau).
%%   \label{dual-scaling-m}
\end{gather}

Modular properties of $\Phi$ are considered after those of the Appell
functions, in~\bref{sec:Phi-S}.

\subsection{Modular transformation properties}
%% The Appell functions $\KK_{\alpha,k}(2\tau,\nu,\mu)$ are
%% $T$-invari\-ant.  
The $S$-transformation of the Appell functions
$\KK_{\alpha,k}(2\tau,\nu,\mu)$ can be derived
from~\eqref{Knew-Stransf-fin}.  We need a version of the
$S$-transformation formula for the functions
$\KKK^a_\alpha(\tau,\nu,\mu)$ introduced in~\bref{sec:S-B}.  The
following lemma plays a crucial role in the calculations in
Sec.~\ref{sec:modular}.\footnote{The lemma also explains the
  usefulness of the $\KKK^a_b$ functions: the sign factor $(-1)^a$ in
  the left-hand side of the formula in the lemma becomes the
  characteristic, reduced to $\{0,1\}$, in the right-hand side.}
\begin{lemma}\label{lemma:basic-S} We have the $S$-transform formula
  \begin{multline*}
    \KKK^a_{\rbar}(-\ffrac{1}{\tau},\ffrac{\nu}{\tau},\ffrac{\mu}{\tau})
%%     =\\
    {}=\tau\,e^{i\pi k\frac{\nu^2 - \mu^2}{2\tau}}\!
    \Bigl(\KK_{\abar,k}(\tau,\nu,\mu) +
    (-1)^{r}e^{i\pi k\mu - i\pi k\frac{\tau}{2}}
    \KK_{\abar,k}(\tau,\nu,\mu\!-\!\tau)\!\Bigr)\\
    +\tau\,e^{i\pi k\frac{\nu^2}{2\tau}}
    \!\!\!\!\sum_{\beta\in\{0,1\}}\sum_{\substack{n=0\\ \nbar=\abar}}^{2k-1}
    (-1)^{r\beta}
    \Phi(2k\tau,k\mu\!-\!n\tau\!-\!\beta k\tau)
    \theta_{n,k}(\tau,\nu).
  \end{multline*}
\end{lemma}
%% \begin{rem}
%% \end{rem}
%% \begin{proof}[Proof of~\bref{lemma:basic-S}]
\begin{proof}
  Several properties of the Appell functions and of the $\Phi$
  function are used here.  First, with~\eqref{theta-rewrite} and
  \eqref{u2},\footnote{And the easily verified property $\K_{k}(\tau,
    \nu + \frac{m}{k}, \mu - \frac{m}{k})= \K_{k}(\tau, \nu, \mu)$,
    $m\in\oZ$.} it readily follows from~\eqref{Knew-Stransf-fin} that
  \begin{multline*}%%%\label{eq:modS-K-1st}
    \KK_{\alpha,k}(-\ffrac{1}{\tau},\ffrac{\nu}{\tau},\ffrac{\mu}{\tau})
    =\ffrac{\tau}{2}\,e^{i\pi k\frac{\nu^2 - \mu^2}{2\tau}}
    \!\!\!\!\sum_{\beta%%%,\gamma
      \in\{0,1\}}\!\!
%%     (-1)^{\alpha \gamma}e^{-i\pi k \gamma^2\frac{\tau}{2} - i\pi k \gamma \mu }
%%     \KK_{\beta + k \gamma, k}(\tau,\nu,\mu + \gamma\tau)
    \KKK^\alpha_\beta(\tau,\nu,\mu)
    \\
    +\ffrac{\tau}{2}\,e^{i\pi k\frac{\nu^2}{2\tau}} \sum_{b=0}^{k-1}
    \sum_{\gamma=0}^1 e^{i\pi\frac{\alpha^2}{2k\tau} -
      i\pi\frac{\alpha b}{k} + i\pi\frac{\mu}{\tau}\alpha}(-1)^{\alpha
      \gamma} \Phi(\ffrac{k\tau}{2},\ffrac{k\mu + \alpha - b\tau}{2})
    \theta_{b + k \gamma,k}(\tau,\nu).
  \end{multline*}
  We next rewrite this for the characteristic $\alpha$ replaced with
  $\alpha+\rbar$ %%%, where $\rbar=[r]_2$,
  and use~\eqref{phi-u} %%%{phi-2k}
  with $u=2$:
  \begin{multline*}%%%\label{eq:modS-K-1st}
    \KK_{\alpha+\rbar,k}(
    -\ffrac{1}{\tau},\ffrac{\nu}{\tau},\ffrac{\mu}{\tau})
    =\ffrac{\tau}{2}\,e^{i\pi k\frac{\nu^2 - \mu^2}{2\tau}}
    \!\!\!\!\sum_{\beta%%%,\gamma
      \in\{0,1\}}\!\!
%%     (-1)^{\alpha \gamma + r\gamma}
%%     e^{-i\pi k \gamma^2\frac{\tau}{2} - i\pi k \gamma \mu }
%%     \KK_{\beta + k \gamma, k}(\tau,\nu,\mu + \gamma\tau)
    \KKK^{\overline{\alpha+r}}_\beta(\tau,\nu,\mu)
    \\[-3pt]
    +\ffrac{\tau}{2}\;e^{i\pi k\frac{\nu^2}{2\tau} +
      i\pi\frac{\alpha\mu}{\tau} + i\pi\frac{\alpha^2}{2k\tau}}
    \!\!\!\!\sum_{\beta,\gamma\in\{0,1\}}\sum_{b=0}^{k-1} (-1)^{\alpha
      \gamma + r\gamma + r\beta}
    e^{-i\pi\frac{\alpha b}{k}}
    \times{}
    \\
    {}\times\Phi(2k\tau,k\mu + \alpha - b\tau - \beta k\tau) \theta_{b
      + k \gamma,k}(\tau,\nu).
  \end{multline*}
  It then follows that
  \begin{multline}\label{K-mod-2}
    \KKK^a_{\alpha+\rbar}(
    -\ffrac{1}{\tau},\ffrac{\nu}{\tau},\ffrac{\mu}{\tau})=
    \\
    = \tau\,e^{i\pi k\frac{\nu^2 - \mu^2}{2\tau}}
    \Bigl(\KK_{\abar,k}(\tau,\nu,\mu) +(-1)^{\alpha+r}e^{-i\pi
      k\frac{\tau}{2} - i\pi k\mu}
    \KK_{\overline{a+k}+k,k}(\tau,\nu,\mu+\tau)\!
    \Bigr)\\
%% +\ffrac{\tau}{2}\; e^{i\pi k\frac{\nu^2}{2\tau} +
%%   i\pi\frac{\alpha\mu}{\tau} + i\pi\frac{\alpha^2}{2k\tau}}
%% \!\!\!\!\sum_{\beta,\gamma\in\{0,1\}}\sum_{b=0}^{k-1} \bigl(1 +
%% (-1)^{a+b+k\gamma}\bigr)
%% e^{-i\pi\frac{\alpha b}{k}}\\*
%% {}\times(-1)^{\alpha \gamma + r\gamma + r\beta}
%% \Phi(2k\tau,k\mu + \alpha - b\tau - \beta k\tau) \theta_{b + k
%%   \gamma,k}(\tau,\nu),
    + \tau\,e^{i\pi k\frac{\nu^2}{2\tau} + i\pi\frac{\alpha\mu}{\tau}
      + i\pi\frac{\alpha^2}{2k\tau}}\,
    X^a_{\alpha,\rbar}(\tau,\nu,\mu),
  \end{multline}
  where we introduce the temporary notation
  \begin{multline*}
    X^a_{\alpha,\rbar}(\tau,\nu,\mu)= \ffrac{1}{2}
    \!\!\sum_{\beta,\gamma\in\{0,1\}}\sum_{b=0}^{k-1} \bigl(1 +
    (-1)^{a+b+k\gamma}\bigr)
    e^{-i\pi\frac{\alpha b}{k}}
    \\
    {}\times(-1)^{\alpha \gamma + r\gamma + r\beta} \Phi(2k\tau,k\mu +
    \alpha - b\tau - \beta k\tau) \theta_{b + k \gamma,k}(\tau,\nu).
  \end{multline*}

  We next observe that by virtue of~\eqref{phi-tau},
  \begin{align*}
    X^a_{\alpha,\rbar}(\tau,\nu,\mu) +{}& (-1)^{r+\alpha}
    \sum_{b=0}^{k-1}\bigl(1+(-1)^{a+b+k}\bigr) e^{-i\pi\frac{\alpha
        b}{k}} e^{-i\pi\frac{(k\mu + \alpha - b\tau)^2}{2k\tau}}
    \theta_{b + k \gamma,k}(\tau,\nu)\\
    ={}& \ffrac{1}{2}\!\sum_{\beta,\gamma\in\{0,1\}}\sum_{b=0}^{k-1}
    \bigl(1 + (-1)^{a+b+k\gamma}\bigr) e^{-i\pi\frac{\alpha b}{k}}
    (-1)^{\alpha \gamma + r\beta}\\*
    &\qquad\qquad\qquad{}\times \Phi(2k\tau,k\mu + \alpha - (b +
    k\gamma)\tau - \beta k\tau)
    \theta_{b + k \gamma,k}(\tau,\nu)\\
    \intertext{and it is easy to see that the equality continues as}
    ={}&\sum_{\beta\in\{0,1\}}\sum_{\substack{n=0\\ \nbar=a}}^{2 k-1}
    e^{-i\pi\frac{\alpha n}{k}} (-1)^{r\beta} \Phi(2k\tau,k\mu +
    \alpha - n\tau - \beta k\tau) \theta_{n,k}(\tau,\nu).
  \end{align*}
  Substituting the $X^a_{\alpha,\rbar}(\tau,\nu,\mu)$ thus expressed
  in~\eqref{K-mod-2}, we also use~\bref{sec:openqp} to replace
  $\KK_{\overline{a+k}+k,k}(\tau,\nu,\mu+\tau)$ with
  $\KK_{\overline{a+k}+k,k}(\tau,\nu,\mu-\tau)$.  After some
  cancellations, this gives
  \begin{multline*}
    \!\!\KKK^a_{\alpha+\rbar}(
    -\ffrac{1}{\tau},\ffrac{\nu}{\tau},\ffrac{\mu}{\tau})
    %% =\\
    =\tau\,e^{i\pi k\frac{\nu^2 - \mu^2\!\!}{2\tau}}\!
    \Bigl(\!\KK_{\abar,k}(\tau,\nu,\mu) +
    (-1)^{\alpha+r}e^{i\pi k\mu - i\pi k\frac{\tau}{2}}
    \KK_{\overline{a+k}+k,k}(\tau,\nu,\mu\!-\!\tau)\!\Bigr)\\
    +\tau\,e^{i\pi k\frac{\nu^2}{2\tau} + i\pi\frac{\alpha\mu}{\tau}
      +i\pi\frac{\alpha^2}{2k\tau}}
    \!\!\!\!\sum_{\beta\in\{0,1\}}\sum_{\substack{n=0\\ \nbar=\abar}}^{2k-1}
    (-1)^{r\beta}
    e^{-i\pi\frac{\alpha n}{k}}
    \Phi(2k\tau,k\mu+\alpha-n\tau-\beta k\tau)
    \theta_{n,k}(\tau,\nu)\\[-4pt]
    +\tau\,e^{i\pi k\frac{\nu^2 - \mu^2}{2\tau}}
    (-1)^{\alpha+r}
    \sum_{\substack{n=k\\ \nbar=\overline{a+k}}}^{2k-1}
    e^{-i\pi\frac{\tau n^2}{2k} + i\pi\mu n}\theta_{n - k,k}(\tau,\nu).
  \end{multline*}

  We finally apply~\eqref{K-char-reduce} to
  $\KK_{\overline{a+k}+k,k}(\tau,\nu,\mu\!-\!\tau)$ in the last
  formula.  Because $\overline{\overline{a+k}+k}=\abar$, we have
  $\overline{a+k}+k=\abar + 2\ell$ with an integer~$\ell$, and
  therefore~\eqref{K-char-reduce} is indeed applicable, with
  $\ell=\half(k+\kbar)-1$ if $\abar=\kbar=1$ and $\ell=\half(k+\kbar)$
  otherwise.  This gives
  \begin{multline*}
    \KKK^a_{\alpha+\rbar}(
    -\ffrac{1}{\tau},\ffrac{\nu}{\tau},\ffrac{\mu}{\tau})
%%     =\\
    {}=\tau\,e^{i\pi k\frac{\nu^2 - \mu^2}{2\tau}}\!
    \Bigl(\KK_{\abar,k}(\tau,\nu,\mu) +
    (-1)^{\alpha+r}e^{i\pi k\mu - i\pi k\frac{\tau}{2}}
    \KK_{\abar,k}(\tau,\nu,\mu\!-\!\tau)\!\Bigr)
    \\
    +\tau\,e^{i\pi k\frac{\nu^2}{2\tau} + i\pi\frac{\alpha\mu}{\tau}
      +i\pi\frac{\alpha^2}{2k\tau}}
    \!\!\!\!\sum_{\beta\in\{0,1\}}\sum_{\substack{n=0\\ \nbar=\abar}}^{2k-1}
    (-1)^{r\beta}
    e^{-i\pi\frac{\alpha n}{k}}
    \Phi(2k\tau,k\mu+\alpha-n\tau-\beta k\tau)
    \theta_{n,k}(\tau,\nu),
  \end{multline*}
  and the identity in the lemma is just the $\alpha=0$ case of this.
\end{proof}

\subsubsection{}\label{sec:Phi-S} We finally quote the
$S$-transformation formula for the $\Phi$ function~\cite{[STT]}:
\begin{equation}\label{Phi-S}
  \Phi(-\ffrac{1}{\tau},\ffrac{\mu}{\tau})
%%   = -i\sqrt{-i\tau}\bigl(
%%   \ex{i\pi\frac{\mu^2}{\tau}}\,
%%   \Phi(\tau,\mu) + 1\bigr)\\
%%   = -i\sqrt{-i\tau}\EX{i\pi\frac{\mu^2}{\tau}}
%%   \Phi(\tau,\mu\!-\!\tau)
  =i\sqrt{-i\tau}\,e^{i\pi\frac{(\mu-1)^2}{\tau}}\Phi(\tau,1-\mu).
\end{equation}

\subsubsection{}\label{sin-Phi-lemma}
We note that all the properties of $\Phi$ in~\bref{sec:Phi} can be
derived directly from the definition as well as from the mere
appearance of $\Phi$ in the $S$-transformation
formula~\eqref{Knew-Stransf-fin} (and from the properties of the
Appell functions).  In particular, \eqref{Phi-S} follows from $S^4=1$
evaluated on the Appell functions (see~\cite{[STT]} for the details).
In application to the higher string functions in this paper, it may be
useful to formulate a more specific identity that ``ensures'' the
relation $S^4=1$ evaluated, e.g., on~$\hB_{m,r}$.  It follows from the
$S$-dual relation to~\eqref{phi-2k},
\begin{equation*}%%%\label{dual-scaling-m}
  \ffrac{1}{2k}\sum_{n=0}^{2k-1}
  \EX{
    2i\pi\frac{n\mu}{\tau}
    \!-\!i\pi\frac{m\,n}{k}
    \!+\!i\pi\frac{n^2}{2k\tau}}
  \Phi(\ffrac{\tau}{2k},\mu\!+\!\ffrac{n}{2k})
  =\Phi(2k\tau,2k\mu\!-\![m]_{2k}\,\tau),
\end{equation*}
and other properties of~$\Phi$ in~\bref{sec:Phi}: for $1\leq m\leq
k-1$ and $k+1\leq m\leq 2k-1$, we have
\begin{equation*}
  2i\sum_{n=1}^{k-1}
  \sin\!\ffrac{\pi m\,n}{k}\;
  \ex{i\pi\frac{n^2}{2k\tau}}
  \Phi(\ffrac{\tau}{2k},\ffrac{n}{2k})
  =-\ffrac{1-(-1)^m\!}{2}\;i\cot\!\ffrac{\pi m}{2 k}
  -\ffrac{i\sqrt{2 k}}{2\sqrt{-i\tau}}
  -2k \Phi(2 k\tau,-m\tau).
\end{equation*}

\section{$\W(k)$ characters in the logarithmic
  $\smash{\protect\hSL2_{k}}$ model~\cite{[S-sl2]}}\label{sec:Wchars}
We here recall the spectral-flow and modular properties of the
characters of the $W$-algebra constructed in~\cite{[S-sl2]}.

\subsection{Spectral-flow properties}\label{sec:sf}
Spectral flow automorphisms act on the character of any
$\hSL2_k$-module~$\mL$ as~\cite{[FSST],[S-sl2]}
\begin{equation}\label{spectral-sl2-general}
  \charSL{\mL}{;\theta}(q,z)=
  q^{\frac{k}{4}\theta^2}\,z^{-\frac{k}{2}\theta}\,
  \charSL{\mL}{}(q,z\,q^{-\theta}).\pagebreak[3]
\end{equation}
The integrable representation characters are well-known to be periodic
under the spectral flow, $\charSL{}{r;1}(q,z)=\charSL{}{p-r}(q,z)$,
and therefore $\charSL{}{r;2}(q,z)=\charSL{}{r}(q,z)$.

Clearly, the rule in~\eqref{spectral-sl2-general} also applies to the
characters of the extended algebra $\W(k)$ of the logarithmically
extended $\hSL2_k$ model.  For $\charW{\pm}{r}$ in~\eqref{chi-plus}
and~\eqref{chi-minus}, calculation then shows that
\begin{align}\label{plus-twisted}
  \charW{+}{r;1}(q,z)&=-\charW{-}{r}(q,z)
  - \omega^-_r(q,z) - \half\charSL{}{p-r}(q,z),
  \\
  \label{minus-twisted}
  \charW{-}{r;1}(q,z)&=-\charW{+}{r}(q,z) - \omega^+_r(q,z)
  \\[-4pt]
  \intertext{for $1\leq r\leq p\!-\!1$ (with $p=k+2$), and,
    similarly,}
  \charW{+}{p;1}(q,z)&=-\charW{-}{p}(q,z) -
  \omega^-_p(q,z),\\
  \charW{-}{p;1}(q,z)&=-\charW{+}{p}(q,z) - \omega^+_p(q,z),
\end{align}
with the $\omega^\pm_r$ given in~\bref{sec:chars}.  On these, the
spectral flow action as in~\eqref{spectral-sl2-general} is in turn
readily evaluated as
\begin{align*}
  \omega^+_{r;1}(q,z)&=-\omega^-_r(q,z)
  -
  \half\charSL{}{p-r}(q,z),\\
  \omega^-_{r;1}(q,z)&=-\omega^+_r(q,z)
  +
  \half\charSL{}{r}(q,z)
\end{align*}
for $1\leq r\leq p\!-\!1$, and
$\omega^{\pm}_{p;1}(q,z)=\omega^{\pm}_{p}(q,z)$.

\subsection{Modular transformation properties}
\label{lemma:mod1}
Under the modular group action, the functions
$\charW{\pm}{r}(\tau,\nu)$ and $\omega^{\pm}_r(\tau,\nu)$, $1\leq
r\leq p$, and $\charSL{}{r}(\tau,\nu)$, $1\leq r\leq p-1$, generate a
representation whose structure can be described as a deformation
%% , via a matrix automorphy factor, 
of the representation
\begin{equation}%%%\label{eq:main-result}
  \Rpi\oplus\oC^2\!\tensor\Rpi
  \,\oplus\,
  \Rmin\oplus\oC^2\!\tensor\Rmin
  \oplus
  \oC^3\!\tensor\Rmin,
\end{equation}
where $\Rmin$ is the $(p\,{-}\,1)$-dimensional $\SLiiZ$ representation
on the integrable $\hSL2_k$ characters, $\Rpi$ is a
$(p\,{+}\,1)$-dimensional representation, $\oC^2$ is the defining
two-dimen\-sional representation, and $\oC^3$ is its symmetrized
square.  The deformation amounts to the occurrence of ``lower''
representation characters in the right-hand side of transformations of
the ``higher''~ones.

\subsubsection{}\label{sl2-modular}The lowest in this sense are the
integrable $\hSL2$-representation characters $\charSL{}{r}$, $1\leq
r\leq p-1$, which span the representation $\Rmin$:
\begin{align}
  \label{min-T-action}
  \charSL{}{r}(\tau+1,\nu)&=
  \lambda_{r,p}\charSL{}{r}(\tau,\nu),
  \qquad\lambda_{r,p}=e^{i\pi(\frac{r^2}{2p}-\frac{1}{4})},\\
  \label{min-S-action}
  \charSL{}{r}(-\ffrac{1}{\tau},\ffrac{\nu}{\tau})
  &=\sqrt{\fffrac{2}{p}}\; e^{i\pi
    k\frac{\nu^2}{2\tau}}\sum_{s=1}^{p-1}
  \sin\!\ffrac{\pi r s}{p}\,\charSL{}{s}(\tau,\nu),
\end{align}

\subsubsection{}\label{minor-modular}Next comes the representation
$\Rpi$ spanned by the linear combinations
\begin{align}
  \minorplus_0(\tau,\nu)&=\omega^-_p(\tau,\nu),\notag
  \\
  \label{minorplus-funcs}
  \minorplus_r(\tau,\nu)&=\omega^+_r(\tau,\nu)+\omega^-_{p-r}(\tau,\nu),
  \quad 1\leq r\leq p\,{-}\,1,
  \\
  \minorplus_p(\tau,\nu)&=\omega^+_p(\tau,\nu),
  \notag
\end{align}
which transform as
\begin{equation}\label{minorplus-S}
  \begin{aligned}
    \minorplus_r(\tau+1,\nu)
    &=\lambda_{r,p}\minorplus_r(\tau,\nu),\\
    \minorplus_r(-\ffrac{1}{\tau},\ffrac{\nu}{\tau})
    &=i\sqrt{\fffrac{2}{p}}\;e^{i\pi k\frac{\nu^2}{2\tau}}
    \Bigl(
    \half\,\minorplus_0(\tau,\nu)
    +\ffrac{(-1)^r}{2}\minorplus_p(\tau,\nu)
    +\sum_{s=1}^{p-1}\cos\!\ffrac{\pi r s}{p}\,\minorplus_s(\tau,\nu)
    \!\!\Bigr)
  \end{aligned}
\end{equation}
for $0\leq r\leq p$.  Next, (a deformation of) the
$\oC^2\!\tensor\Rmin$ representation is spanned by the linear
combinations
\begin{equation}\label{skewminor-funcs}
  \begin{split}
    \skewminor_r(\tau,\nu) &=
    (p-r)\omega^+_r(\tau,\nu) - r\omega^-_{p-r}(\tau,\nu),
    \\
    \tauskewminor_r(\tau,\nu) &= \tau\skewminor_r(\tau,\nu),
  \end{split}
  \qquad 1\leq r\leq p\!-\!1,
\end{equation}
which transform as
\begin{align}%%%\label{skewminor-T}
  \skewminor_r(\tau+1,\nu)&=\lambda_{r,p} \skewminor_r(\tau,\nu),
  \quad \tauskewminor_r(\tau+1,\nu)=\lambda_{r,p}
  \bigl(\tauskewminor_r(\tau,\nu) +\skewminor_r(\tau,\nu)\bigr),
  \notag
  \\
  \label{skewminor-S}
  \skewminor_r(-\ffrac{1}{\tau},\ffrac{\nu}{\tau})
  &=\sqrt{\fffrac{2}{p}}\, e^{i\pi k\frac{\nu^2}{2\tau}}
  \sum_{s=1}^{p-1}\sin\!\ffrac{\pi r s}{p}
  \Bigl(\tauskewminor_s(\tau,\nu)
  - \ffrac{p\nu}{2}\charSL{}{s}(\tau,\nu)\Bigr),
  \\
  \tauskewminor_r(-\ffrac{1}{\tau},\ffrac{\nu}{\tau})
  &=\sqrt{\fffrac{2}{p}}\, e^{i\pi k\frac{\nu^2}{2\tau}}
  \sum_{s=1}^{p-1}\sin\!\ffrac{\pi r s}{p}
  \Bigl(-\skewminor_s(\tau,\nu) +
  \ffrac{p\nu}{2\tau}\charSL{}{s}(\tau,\nu)\Bigr)
  \notag
\end{align}
(the deformation is due to $\nu$ times the integrable-representation
characters occurring in the right-hand
side).

\subsubsection{}\label{chi-modular}Further, the linear combinations of
the $\W(k)$-characters
\begin{align*}
  \rho_0(\tau,\nu)&=\charW{-}{p}(\tau,\nu),\\[-3pt]
  \rho_r(\tau,\nu)&=\charW{+}{r}(\tau,\nu) + \charW{-}{p-r}(\tau,\nu)
  + \ffrac{r}{2p}\charSL{}{r}(\tau,\nu),\quad
  1\leq r\leq p\!-\!1,\\[-3pt]
  \rho_p(\tau,\nu)&=\charW{+}{p}(\tau,\nu)
\end{align*}
transform as
\begin{align*}
  \rho_r(\tau+1,\nu)&=\lambda_{r,p}\rho_r(\tau,\nu),
  \\
  \rho_r(-\ffrac{1}{\tau},\ffrac{\nu}{\tau}) &=i\sqrt{\fffrac{2}{p}}\,
  e^{i\pi k\frac{\nu^2}{2\tau}}
  \biggl(\ffrac{1}{2}(\tau\rho_0(\tau,\nu) +
  \nu\minorplus_0(\tau,\nu)) +\ffrac{(-1)^r}{2}(\tau\rho_p(\tau,\nu) +
  \nu\minorplus_p(\tau,\nu))\\
  &\qquad\qquad\qquad{}+\sum_{s=1}^{p-1} \cos\!\ffrac{\pi r s}{p}\,
  \bigl(\tau\rho_s(\tau,\nu) + \nu
  \minorplus_s(\tau,\nu)\bigr)\!\!\biggr).
\end{align*}
Here, $\tau\rho_r(\tau,\nu)$ are to be regarded as new functions, with
the modular transformations for them to be (easily) obtained from the
above formulas (for example, $\tau\rho_r\mapsto
\lambda_{r,p}\tau\rho_r + \lambda_{r,p}\rho_r$ under $\tau\mapsto
\tau+1$; we do not introduce a special notation for $\tau\rho_r$).
Modulo the $\nu$-terms in the right-hand sides, $(\rho_r,\tau\rho_r)$
then span the $\SLiiZ$ representation $\oC^2\!\tensor\Rpi$.

Finally, the linear combinations of the characters
\begin{equation*}
  \skewmajor_r(\tau,\nu)
  =(p-r)\charW{+}{r}(\tau,\nu) - r\charW{-}{p-r}(\tau,\nu)
  -\Bigl(\ffrac{r^2}{4p}
  +\ffrac{1}{8i\pi\tau}
  \Bigr)\charSL{}{r}(\tau,\nu),
  \quad 1\leq r\leq p\!-\!1,
\end{equation*}
transform as
\begin{align*}
  \skewmajor_r(\tau+1,\nu)&=\lambda_{r,p}\skewmajor_r(\tau,\nu),
  \\
  \skewmajor_r(-\ffrac{1}{\tau},\ffrac{\nu}{\tau})
  &=\sqrt{\fffrac{2}{p}}\, e^{i\pi k\frac{\nu^2}{2\tau}}
  \sum_{s=1}^{p-1}\sin\!\ffrac{\pi r s}{p}\,
  \Bigl(\tau^2\skewmajor_s(\tau,\nu) + \nu\tau\skewminor_s(\tau,\nu) -
  \ffrac{p\nu^2}{4} \charSL{}{s}(\tau,\nu)\!\Bigr).
\end{align*}
Here, too, $(\skewmajor_r,\tau\skewmajor_r,\tau^2\skewmajor_r)$ form
the triplet $\oC^3\!\tensor\Rmin$ modulo the explicitly
$\nu$-dependent terms.

\section{$\hA\hB\hC$ identities}\label{app:ABC}
We here derive the ``open periodicity'' and some other symmetries of
the higher string functions.

\begin{Lemma}\label{app:theF}%
  \addcontentsline{toc}{subsection}{\thesubsection.  \ \ The
    $\hA\hB\hC$ lemma} For
  \begin{align}\label{the-C}
    C_{n,r}(q)&=\sum_{a\in\oZ}\sum_{j\geq 1}(-1)^{j+1}
    \Bigl(q^{\half j(j-n) + r a + p a^2 + \half j(2a p + r)}
    -(r\mapsto-r)\!\Bigr),\\
    \label{the-B}
    B_{n,r}(q)&=\sum_{a\in\oZ}\sum_{j\geq 1}(-1)^{j+1}
    a\Bigl(q^{\half j(j-n) + r a + p a^2 + \half j(2a p + r)}
    -(r\mapsto-r)\!\Bigr),\\
    \label{the-A}
    A_{n,r}(q)
    &= \sum_{a\in\oZ}\sum_{j\geq 1}(-1)^{j+1}
    a^2 \Bigl(q^{\half j(j-n) + r a + p a^2 + \half j(2a p + r)}
    -(r\mapsto-r)\!\Bigr),
  \end{align}
  we have the ``open quasiperiodicity'' formulas
  \begin{multline*}
    C_{n+2k\ell,r}(q)=q^{k\ell^2 + n\ell}\,C_{n,r}(q),\\[-1.2\baselineskip]
  \end{multline*}
  \begin{multline*}
    B_{n+2k\ell,r}(q)=
    q^{k\ell^2 + n\ell}\,B_{n,r}(q)
    +\ell\,q^{k\ell^2 + n\ell}\,C_{n,r}(q) +{}
    \\
    {}+
    \begin{cases}\displaystyle
%%       +
      -
      q^{k\ell^2 + n\ell}
      \sum_{j=1}^{2\ell}(-1)^{j}q^{-\frac{k}{4}j^2 - \frac{n}{2}j
        - \frac{r^2}{4p}}\,
%%       \Psi_{r}^j(q)
      \eta(q)\psi^{\jbar}_{r}(q)
      ,& \ell\geq 1,
      \\
      \displaystyle
%%       -
      q^{k\ell^2 + n\ell}
      \!\!\sum_{j=2\ell+1}^{0}(-1)^{j}q^{-\frac{k}{4}j^2 - \frac{n}{2}j
        - \frac{r^2}{4p}}
%%       \,\Psi_{r}^j(q)
      \eta(q)q^{}\,\psi^{\jbar}_{r}(q)
      ,& \ell\leq -1,
    \end{cases}
  \end{multline*}
  and
  \begin{multline*}
    A_{n+2k\ell,r}(q)=
    q^{k\ell^2 + n\ell}\,A_{n,r}(q)
    + 2\ell\,q^{k\ell^2 + n\ell}\,B_{n,r}(q)
    + \ell^2\,q^{k\ell^2 + n\ell}\,C_{n,r}(q)+{}
    \\
    {}+
    \begin{cases}\displaystyle
%%       +
      -
      q^{k\ell^2 + n\ell}\sum_{j=1}^{2\ell}(2\ell-j)(-1)^{j}
      q^{-\frac{k}{4}j^2 - \frac{n}{2}j
        - \frac{r^2}{4p}}
%%       \,\Psi_{r}^j(q)
      \eta(q)\,\psi^{\jbar}_{r}(q)
      ,&\ell\geq 1,
      \\
      \displaystyle
%%       -
      q^{k\ell^2 + n\ell}\!\!\sum_{j=2\ell+1}^{0}\!(2\ell-j)(-1)^{j}
      q^{-\frac{k}{4}j^2 - \frac{n}{2}j
        - \frac{r^2}{4p}}
%%       \,\Psi_{r}^j(q)
      \eta(q)\,\psi^{\jbar}_{r}(q)
      ,&\ell\leq -1,
    \end{cases}
  \end{multline*}
  where
  \begin{equation}\label{the-Psi}
    \psi^{\jbar}_r(q)=
    \begin{cases}
      \psi^+_q(q),& j\text{ even},\\
      \psi^-_q(q),& j\text{ odd}.
    \end{cases}
  \end{equation}
\end{Lemma}

\noindent
Definition~\eqref{the-Psi} is an excusable abuse of notation.  The
formula for $C_{n+2k\ell,r}(q)$ is of course the classic
string-function ``quasiperiodicity.''

\begin{proof}
  The properties claimed in the lemma are particular cases of a
  general formula derived as follows.  For a (polynomial) function $f$
  defined on $\oZ$, we set
  \begin{equation*}
    F_{n,r}(q)=\sum_{a\in\oZ}\sum_{j\geq1}(-1)^{j+1}f(a)
    \Bigl(q^{\half j(j-n) + r a + p a^2 +\half j(2 a p + r)}
    -(r\mapsto-r)\!\Bigr)
  \end{equation*}
  and then calculate $F_{n + 2 k\ell,r}$ for any $\ell\in\oZ$ by
  shifting the summation variables as $a\mapsto a+\ell$ and $j\mapsto
  j-2\ell$.  An elementary calculation then gives
  \begin{multline*}
    F_{n + 2 k\ell,r}(q)=q^{k\ell^2 + n\ell}
    \sum_{a\in\oZ}\sum_{j\geq 2\ell+1}\!\!(-1)^{j+1}f(a+\ell)
    \Bigl(q^{\half j(j-n) + r a + p a^2 + \half j(2 a p + r)}
    -(r\mapsto-r)\!\Bigr)\\
    \shoveleft{{}
      =q^{k\ell^2 + n\ell} F_{n,r}}\\*
    {}+
    q^{k\ell^2 + n\ell}
    \sum_{a\in\oZ}\sum_{j\geq1}(-1)^{j+1}
    \bigl(f(a+\ell)-f(a)\bigr)
    \Bigl(q^{\half j(j-n) + r a + p a^2 + \half j(2 a p + r)}
    -(r\mapsto-r)\!\Bigr)\\
    +q^{k\ell^2 + n\ell}
    \sum_{a\in\oZ}
    \begin{bmatrix}
      -\!\!\sum\limits_{j=1}^{2\ell}\\
      \sum\limits_{j=2\ell+1}^{0}
    \end{bmatrix}
    (-1)^{j+1}
    f(a+\ell)
    \Bigl(q^{\half j(j-n) + r a + p a^2 + \half j(2 a p + r)}
    -(r\mapsto-r)\!\Bigr),
  \end{multline*}
  where $-\!\!\sum\limits_{j=1}^{2\ell}$ is to be taken for $\ell>0$
  and $\sum\limits_{j=2\ell+1}^{0}$ for $\ell<0$.  In either of these
  \textit{finite} sums, we can change the order of summation and then
  shift the $a$ summation variable to obtain
  \begin{multline*}
    F_{n + 2 k\ell,r}(q)=q^{k\ell^2 + n\ell}\Bigl(F_{n,r}+\\
    {}+
    \sum_{a\in\oZ}\sum_{j\geq1}(-1)^{j+1}
    \bigl(f(a+\ell)-f(a)\bigr)
    \Bigl(q^{\half j(j-n) + r a + p a^2 + \half j(2 a p + r)}
    -(r\mapsto-r)\!\Bigr)\\
    + \begin{bmatrix}
      -\!\!\sum\limits_{j=1}^{2\ell}\\
      \sum\limits_{j=2\ell+1}^{0}
    \end{bmatrix}
    (-1)^{j+1}q^{-\frac{k}{4}j^2 - \frac{n}{2}j}
    \!\!\sum_{a\in\oZ+\frac{j}{2}}\!\!
    \bigl(f(\ell - \ffrac{j}{2} + a) - f(\ell - \ffrac{j}{2} - a)\bigr)
    q^{p a^2 + r a}\Bigr).
  \end{multline*}
  For $f(a)=1$, $a$, and $a^2$, the respective ``$F$''-functions are
  $C_{n,r}(q)$, $B_{n,r}(q)$, and $A_{n,r}(q)$, with the results
  stated in the lemma.
\end{proof}

\subsubsection{}\label{hBAshift}
For $\hC_{n,r}(q)$, $\hB_{n,r}(q)$, and $\hA_{n,r}(q)$ expressed as
in~\eqref{ABC-help}, the formulas in~\bref{app:theF} are restated as
follows: first, $\hC_{n+2k\ell,r}(q)=\hC_{n,r}(q)$, and then
\begin{align*}
  \hB_{n+2k\ell,r}(q)&=
  \hB_{n,r}(q)
  +
  \begin{cases}\displaystyle
%%     +
    -
    \sum_{j=1}^{2\ell}(-1)^{j}q^{-\frac{k}{4}(j+\frac{n}{k})^2}
    \ffrac{\psi_{r}^{\jbar}(q)}{\eta(q)}
    ,& \ell\geq 1,
    \\
    \displaystyle
%%     -
    \!\!\sum_{j=2\ell+1}^{0}(-1)^{j}q^{-\frac{k}{4}(j+\frac{n}{k})^2}
    \ffrac{\psi_{r}^{\jbar}(q)}{\eta(q)},& \ell\leq -1,
  \end{cases}
  \\
  \intertext{and}
  \hA_{n+2k\ell,r}(q)&=
  \hA_{n,r}(q)+
  \begin{cases}\displaystyle
%%     -
    \sum_{j=1}^{2\ell}(j+\ffrac{n}{k})(-1)^{j}
    q^{-\frac{k}{4}(j+\frac{n}{k})^2}
    \ffrac{\psi_{r}^{\jbar}(q)}{\eta(q)},&\ell\geq 1,
    \\
    \displaystyle
%%     +    
    -
    \!\!\sum_{j=2\ell+1}^{0}\!(j+\ffrac{n}{k})(-1)^{j}
    q^{-\frac{k}{4}(j+\frac{n}{k})^2}
    \ffrac{\psi_{r}^{\jbar}(q)}{\eta(q)},&\ell\leq -1
  \end{cases}
\end{align*}
(we recall that $\psi_{r}^{\jbar}(q)$ is defined in~\eqref{the-Psi}).

%% \subsection{Some other $\hA\hB\hC$ symmetries}
\begin{Lemma}\label{sec:ABC-sym}%
  \addcontentsline{toc}{subsection}{\thesubsection.  \ \ Some other
    $\hA\hB\hC$ symmetries}%
  Relations~\eqref{B-minus}, \eqref{A-minus}, \eqref{sfB-properties},
  and~\eqref{sfA-properties} hold.
\end{Lemma}
\begin{proof}
  The reflection symmetries, Eqs.~\eqref{B-minus} and~\eqref{A-minus},
  are shown by elementary manipulations with the same $F_{n,r}(q)$ as
  in~\bref{app:theF}, which yield
  \begin{multline*}
    F_{-n,r}(q)=\sum_{a\in\oZ}\sum_{j\geq1}(-1)^{j+1}f(-a)
    \Bigl(q^{\half j(j-n) + r a + p a^2 +\half j(2 a p + r)}
    -(r\mapsto-r)\!\Bigr)
    \\[-6pt]
    {}+\sum_{a\in\oZ}\bigl(f(a)-f(-a)\bigr)q^{r a + p a^2}
  \end{multline*}
  (identity~\eqref{zero-dentity} was used here in particular).  For
  $f(a)=1$, we recover the well-known symmetry
  $C_{-n,r}(q)=C_{n,r}(q)$ and, evidently,
  $\hC_{-n,r}(q)=\hC_{n,r}(q)$ for $\hC_{n,r}(q)$ defined
  in~\eqref{sec:results}; Eqs.~\eqref{B-minus} and~\eqref{A-minus}
  also follow immediately.

  Next, the first line in~\eqref{sfB-properties} follows
  from~\bref{hBAshift} and~\eqref{phi-tau}, and the second line
  from~\eqref{B-minus}, \eqref{Phi-reflect}, and~\eqref{phi-tau}.
  Similarly, the first line in~\eqref{sfA-properties} follows
  from~\bref{hBAshift} and~\eqref{phi-tau}, and the second line
  from~\eqref{A-minus}, \eqref{phi-tau}, and the identity
  \begin{equation*}
    \Phi'(\tau,-\mu) =
    \Phi'(\tau,\mu) - \ffrac{\mu}{\tau}\;e^{-i\pi\frac{\mu^2}{\tau}}
  \end{equation*}
  obtained by differentiating~\eqref{Phi-reflect}
  (see~\eqref{Phi-prime}).
\end{proof}

\subsubsection{}With relations~\eqref{B-minus} thus established, it
readily follows from~\bref{hBAshift} that modulo
$\oC[\psi^{\pm}_s(q)q^{-\frac{n^2}{4 k}}/\eta(q)]$, the independent
$\hB_{n,r}(q)$ are
%% \begin{equation}\label{eq:B-independent}
  $\hB_{m,r}(q)$,  $1\leq m\leq k-1$, $1\leq r\leq p$.
%%   \pagebreak[3]
%% \end{equation}
In particular, it is easy to see that\pagebreak[3]
\begin{equation}\label{B-eval}
  \hB_{-k,r}(q)=-\ffrac{\psi^-_r(q)}{2\eta(q)},\qquad
  \hB_{0,r}(q)=-\ffrac{\psi^+_r(q)}{2\eta(q)},
\end{equation}
and so on for $\hB_{k\ell,r}(q)$ in accordance with~\bref{hBAshift}.

\subsubsection{}Finally, it is also obvious from the definitions
in~\bref{sec:results} that
%% \begin{equation*}
  $\hC_{n,0}(q)=\hB_{n,0}(q)=\hA_{n,0}(q)=0$.
%% \end{equation*}
In view of the symmetry
\begin{equation}\label{C-mirror}
  \hC_{n+k,p-r}(q) = \hC_{n,r}(q),
\end{equation}
this also implies that~$\hC_{n,p}(q)=0$.

%% Next,
%% \begin{equation*}
%%   \hB_{n+k,p-r}(q)=\hB_{n,r}(q)
%%   -\ffrac{\psi^+_{p-r}(q)}{\eta(q)}\,q^{-\frac{(n+k)^2}{4k}}
%%   +\hX_{n,r}(q),
%% \end{equation*}
%% where
%% \begin{equation*}
%%   \hX_{n,r}(q)=\half\,\ffrac{q^{-\frac{n^2}{4k}}}{\eta(q)^2}
%%   \sum_{a\in\oZ}\sum_{j\geq 1}(-1)^{j+1}
%%   \Bigl(q^{\half j(j-n) + \frac{(2 a p + r)^2}{4p} + \half j(2a p + r)}
%%   +(r\mapsto-r)\!\Bigr).
%% \end{equation*}
%% It also
%% follows that
%% \begin{multline*}
%%   \sfBfirst_{m+k,p-r}(\tau)
%%   =
%%   \sfBfirst_{m,r}(\tau)
%%   + \hX_{m,r}(\tau)
%%   +\ffrac{\theta_{r,p}(\tau)}{\eta(\tau)^2}\,\Phi(2k\tau, m\tau)
%%   -\ffrac{\theta_{p-r,p}(\tau)}{\eta(\tau)^2}\,\Phi(2k\tau, (m-k)\tau).
%% \end{multline*}

\end{document}